\documentclass[11pt]{amsart}
\usepackage[T1]{fontenc}
\usepackage{amssymb,amsmath,amsthm}

\newtheorem{theorem}{Theorem}[subsection]

\newtheorem{claim}[theorem]{Claim}

\newtheorem{corollary}[theorem]{Corollary}

\newtheorem{definition}{Definition}[section]

\newtheorem{lemma}[theorem]{Lemma}

\newtheorem{remark}[theorem]{Remark}

\newenvironment{corollarytotheproof}[1][Corollary \emph{(to the 
proof)}.]{\begin{trivlist}
\item[\hskip \labelsep {\bfseries #1}]}{\end{trivlist}}

\begin{document}

\def\Ind#1#2{#1\setbox0=\hbox{$#1x$}\kern\wd0\hbox to 0pt{\hss$#1\mid$\hss}
\lower.9\ht0\hbox to 0pt{\hss$#1\smile$\hss}\kern\wd0}
\def\ind{\mathop{\mathpalette\Ind{}}}
\def\Notind#1#2{#1\setbox0=\hbox{$#1x$}\kern\wd0\hbox to 0pt{\mathchardef
\nn=12854\hss$#1\nn$\kern1.4\wd0\hss}\hbox to
0pt{\hss$#1\mid$\hss}\lower.9\ht0 \hbox to
0pt{\hss$#1\smile$\hss}\kern\wd0}
\def\nind{\mathop{\mathpalette\Notind{}}}

\def\thind{\mathop{\mathpalette\Ind{}}^{\text{\th}} }
\def\nthind{\mathop{\mathpalette\Notind{}}^{\text{\th}} }
\def\uth{\text{U}^{\text{\th}} }

\def\spec{\mathop{Spec}}


\title{Properties and Consequences of Thorn-Independence}
\author{Alf Angel Onshuus}


\maketitle

\begin{abstract}
In the last couple of decades, independence relations have become one
of the central parts of model theory. Areas such as classification theory,
stability theory and simplicity theory use the notion of
``forking independence''
defined by Shelah  as a main tool for
their development. However, as powerful and useful as forking-independence
is, it lacks generality. O-minimal structures give an example of a whole
class of structures for which forking independence does not work. This, we
believe, is the main cause why even though the results in o-minimal theory
seem very similar to those in the theory of strongly minimal
sets (strongly minimal
sets are in some ways the smallest structures for which forking defines an
independence relation) the proofs of analogous results have been very
different in both areas.

We begin by developing a new notion of independence
(\th-independence, read ``thorn''-independence) that arises from a
family of ranks suggested by Scanlon (\th-ranks). We prove that in
a large class of theories (it includes simple theories and
o-minimal theories) this notion has many of the properties needed
for an adequate geometric structure.

Finally, we analyze the behavior of \th-forking in some theories
where an independence relation had already been studied by other
authors. We prove that \th-independence agrees with the usual
independence notions in stable, supersimple and o-minimal
theories. Furthermore, we give some evidence that the equivalence
between forking and \th-forking in simple theories might be
closely related to one of the main open conjectures in simplicity
theory, the stable forking conjecture. In particular, we prove
that in any simple theory where the stable forking conjecture
holds, \th-independence and forking independence agree.
\end{abstract}






\section{Introduction}

\subsection{Overview}

In the last decades independence notions have become one of the central
ideas in model theory. Especially in the last three decades independence notions and ranks have been the key factor first in Shelah's classification theory and later in understanding
the geometry in the models of theories such as stable, o-minimal
and in the last decade simple theories. However, even when many of
the results on the topology of o-minimal theories were similar to some of those in geometric stability theory,
they were developed almost independently from each other mainly
because of the difference between the definitions of the respective
independence notions.

In sections 2 through 4 we develop a new notion of independence
(\th-independence --read ``thorn'' independence--) and a rank that
is associated to this independence notion (\th-rank). We prove
that in a large class of theories which we shall call ``rosy''
this notion defines a geometric independence relation. This class
of theories includes simple and o-minimal theories. It also
includes theories for which there was no previously known
independence relation; for example models of the theory defined by
Casanovas and Wagner in ~\cite{Casanovas/Wagner}\footnote{ In this
paper, Casanovas and Wagner construct a theory which is not simple
and provides the first example of a theory without the strict
order property which does not eliminate hyperimaginaries} are
rosy.

In section 5 we study the relation between the classical
independence notions which we had in simple and o-minimal theories
and \th-independence. We prove that \th-independence agrees with
the classical independence notion in o-minimal, stable and all
known cases of simple theories. This provides a unified approach
in two areas where, until now, the methods and proofs have been
different (even though, as mentioned above, they did go in the
same general direction and similar results were obtained). It also
gives an alternative (and, in our point of view simplified)
definition of forking in stable theories which might give further
insight in stability theory.

\subsection{Notation and Conventions}

We assume the reader is familiar with the terminology and the
basic results of model theory and, more specifically, stability
and simplicity theory.

As it is common in stability theory, given a complete theory $T$
we will fix a universe $\mathcal{C}$ called a ``monster'' model of
$T$: we choose some saturated model $\mathcal{C}$ of cardinality
$\kappa$ and assume all sets, types and models we talk about have
cardinality less than $\kappa$ and live inside $\mathcal{C}$. In
particular, by \emph{models of $T$} we mean an elementary submodel
of $\mathcal{C}$. Any automorphism will be understood to be a
$\mathcal{C}$-automorphism. Following \cite{hodges}, an
\emph{equivalence formula} of $\mathcal{C}$ is a formula
$\phi(x,y)$ in the language of $T$ that defines an equivalence
relation in $\mathcal{C}$. Unless otherwise specified, we will
work inside $\mathcal{C}^{eq}$ in the sense of \cite{Shelah}. By
convention lower case letters $a,b,c,d$ will in general represent
tuples (of imaginaries) and upper case letters will represent
sets. Greek letters such as $\delta, \sigma, \psi, \phi$ will be
used for formulas.

\section{\th-Forking}

\subsection{Definitions}

We will start by defining the notions that we will work with
throughout this paper.

\begin{definition}

A formula $\delta (x,a)$ \emph{strongly divides over $A$} if
$tp(a/A)$ is non-algebraic and $\{ \delta (x,a') \}_{a'\models
tp(a/A)}$ is $k$-inconsistent for some $k\in \mathbb{N}$.

We will say that $\delta (x,a)$ \emph{\th-divides over $A$} if we
can find some tuple $c$ such that $\delta (x,a)$ strongly divides
over $Ac$.

Finally, a formula \emph{\th-forks} over $A$ if it implies a
(finite) disjunction of formulae which \th -divide over $A$.
\end{definition}

\begin{remark}
Suppose a formula $\delta(x,a)$ strongly divides over some set
$C$; let $p(x,C):=tp(a/C)$. By definition there is some $k\in
\mathbb{N}$ such that

\noindent $\{ \delta (x,a') \}_{a'\models tp(a/C)}$ is
$k$-inconsistent.

Another way of saying this is that for all $x_1,x_2, \dots ,x_k$
\[ p(x_1,C)\cup p(x_2,C) \cup \dots \cup p(x_k,C) \models \neg \left( \bigwedge_{i=1}^n \phi(x,a_i) \right).\]

By compactness, there is some formula $\theta (y,c) \in tp(a/C)$
such that

\noindent $\{ \delta(x,a') \}_{a'\models \theta(y,c)}$ is
$k$-inconsistent.
\end{remark}

We will say the type $p(x)$ \th-divides over $A$ if there is a
formula in $p(x)$ which \th-divides over $A$; similarly for
\th-forking. We say that $a$ is \emph{\th-independent} of $b$ over
$A$, denoted $a\thind_A b$, if $tp\left( a/Ab\right) $ does not
\th -fork over $A$.

Note that even though \th-dividing and dividing have similar
definitions, the fact that we ask for $k$-incosnsitency for a set
of formulas for which the parameters vary in a definable class as
opposed to an indiscernible sequence is a significant
modification. In particular, as we will prove later, many theories
with the strict order property will behave nicely under this new
definition.

\subsubsection{First Results}

\begin{lemma}
\label{introduction} Let $A,B,C$ be subsets (of $\mathcal{C}$)
such that $A\subseteq B\subseteq C$ and let $a,b$ be tuples. Then
\begin{enumerate}
\item \emph{Extension:} Given a type $p$ over $B$ which does not
\th-fork over $A$, we can extend $p$ to a type $q$ over $C$ which
does not \th-fork over $A$.

\item \emph{Monotonicity:} If a formula $\delta(x,b)$ \th-forks
over $B$, then it \th-forks over $A$.

\item \emph{Partial right transitivity:} If $tp\left( a/C\right) $
does not \th-fork (divide) over $A$ then $tp\left( a/C\right) $
does not \th-fork (divide) over $B$ and $tp\left( a/B\right) $
does not \th-fork (divide) over $A.$

\item If $b\thind _{A} a$ and $tp\left( a/A\right) $ is
non-algebraic, then $tp\left( a/Ab\right) $ is not algebraic.

\item If $b\thind _{A} a$ and $tp\left( b/A\right) $ is
non-algebraic, then $tp\left( b/Aa\right) $ is not algebraic.

\item If $p(x,b)$ is a type over $Ab$ which does not \th-fork over
$A$ and $\theta (x,y)$ is a formula such that $p(x,b) \models
\exists y \theta (x,y)$. Then $p(x,b)\cup \{ \theta (x,y) \} $
does not strongly divide over $A$.

\item \emph{Base Extension:} If $a\thind _A c$, then for any tuple
$d$ there is some $d'\models tp(d/Ac)$ such that $a\thind _{Ad'}
c$.

\end{enumerate}
\end{lemma}

\begin{proof}.

\noindent \emph{1)} Let $\Gamma (x):=\{ \psi(x) | \psi(x)\in
\mathcal{L}(x) \text{ and $\neg \psi(x)$ \th-divides over $A$}\}$.
As in the proof of extension for forking, we first need to prove
that the type $p(x)\cup \Gamma(x)$ is consistent. Suppose this was
not the case. By compactness there would be a finite subset
$\psi_1(x), \dots \psi_n(x)$ in $\Gamma(x)$ such that $p(x)\models
\vee_{i=1}^n \neg \psi_i(x)$. By definition $p(x)$ would \th-fork
over $A$. Let $q'$ be a complete extension of $p(x)\cup \Gamma(x)$
to some set containing $C$, and let $q$ be the restriction of $q'$
to $C$. By definition $q$ does not \th-fork over $A$.

\noindent \emph{2)} In the definition of \th-dividing we are
allowed to add parameters to get strong division, so if any
formula \th-divides over an extension of $A$ then it \th-divides
over $A$. The result for \th-forking follows.

\noindent \emph{3)} The first implication follows from
monotonicity and the other follows from the fact that we are
considering fewer formulas.

\noindent \emph{4)} Let $p(x,b)=tp(a/Ab)$. Assuming $p(x,b)$ does
not strongly divide and using non-algebraicity of $tp(b/A)$, we
can find distinct realizations $b_1, b_2, \dots ,b_n, \dots $ of
$tp(b/A)$ such that $\bigcup_i p(x,b_i)$ is consistent, realized
by some $a'$. However, as $tp(b_1/A)=tp(b/A)$ we can assume
$a'b_1=ab$. However, $p(x,b)$ is a complete type, so $b_i\models
tp(b/Aa)$ showing that $tp(b/Aa)$ is non algebraic.

\noindent \emph{5)} We will prove the contrapositive of the
statement. Let $tp(b/Aa)$ be algebraic and let $b=b_1, b_2, \dots
b_n$ be all the elements that satisfy $tp(b/Aa)$. Then
\[ tp(b/Aa)\models (x=b_1)\vee (x=b_2)
\vee \dots \vee (x=b_n).\]  If $tp(b/A)$ is non-algebraic, then
$x=b_i$ strongly divides (and therefore \th-divides) over $A$ for
all $i\leq n$; by definition $tp(b/Aa)$ \th-forks over $A$ and
$b\nthind _A a$.

\noindent \emph{6)} Let us suppose there is some $\delta (x,b) \in
p(x,b)$ such that $\delta (x,b) \wedge \theta (x,y)$ strongly
divides over $A$ and suppose $\delta (x,b) \wedge \exists y \theta
(x,y)$ does not strongly divide over $A$. By definition (and
compactness) there are some distinct $b_1,\dots , b_n ,\dots $
such that   $b_i\models tp(b/A)$ and $\bigwedge_{i\in \mathbb{N}}
\delta (x,b_i)\wedge \exists y \theta (x,y)$ is consistent,
satisfied by some element $a$. Let $c\models \theta(a,y)$. Then
$(a,c)\models \delta (x,b_i) \wedge \theta (x,y)$ for all $i$,
contradicting our assumptions.

\noindent \emph{7)} By extension we can find some $a'\models
tp(a/Ac)$ such that $a'\thind _A cd$. Let $d'$ be the image of an
automorphism which sends $a'$ to $a$ fixing $Ac$, so that $a\thind
_A cd'$. By partial right transitivity, $a\thind _{Ad'} c$.
\end{proof}


\begin{remark}
Condition \emph{5} above does not hold for \th-dividing in place
of \th-forking, even in the theory in the language of equality
that states that there are infinitely many elements. In fact, if
$A$ is the unordered pair $\{a,b\}$ where $a\leq b$, then
$tp(a/A)$ is algebraic, but $x\in A$ does not \th-divide over the
empty set.
\end{remark}
\begin{proof} We will prove the remark by contradiction.
Let $B$ be any set such that $tp(A/B)$ is non-algebraic and $x\in
A$ strongly divides over $B$. If $tp(a/Bb)$ was non-algebraic,
then $\{ x\in A' \}_{A'\models tp(A/B), b\in A'}$ would be an
infinite consistent set of formulas which would contradict strong
division. This means that $b\in {\rm acl}(Ba)$; by symmetry, $a\in
{\rm acl}(Bb)$. In this particular theory we may assume that $B$
is a subset of $M$ (as opposed to $M^{\rm eq})$\footnote{This is
because the theory of equality has ``elimination of imaginaries'',
a property we will talk about further on}. But $tp(A/B)$ is
non-algebraic, so $tp(b/B)$ and $tp(a/B)$ must be both be
non-algebraic. By Steinitz exchange principle (see for example
\cite{Pillay}, definition 2.1.1), there must be some element in
$B$ algebraic over $a$. However, ${\rm acl}(a)=\{a\}$, so $a\in B$
which contradicts $a$ being non-algebraic over $B$.
\end{proof}


This shows that, unlike the case of simple theories, even when we
are only considering the language of equality there is a
difference between \th-dividing and \th-forking.

The following results go in the same direction as the ones above
(showing properties that \th-forking has in a general theory), but
the proofs are a little more elaborated.

\begin{theorem}
Given a tuple $a$ and sets $A\subset B$, then $tp(a/B)$ \th-forks
over $A$ if and only if $tp(a/acl(B))$ \th-forks over $A$. In
other words, $a\thind _A B$ if and only if $a\thind _A acl(B)$.
\end{theorem}

\begin{proof}
Monotonicity implies that whenever $tp(a/B)$ \th-forks over $A$,
$tp(a/acl(B))$
  \th-forks over $A$. For the other direction, suppose $\psi(x,c)\in tp(a/acl(B))$ is such that $\psi(x,c)$ \th-forks over $A$. Let $\phi_i(x,d_i)$ be such that \[ \psi(x,c)\Rightarrow \phi_i(x,d_1) \vee \phi(x,d_2) \vee \dots \vee \phi(x,d_m)\] and for any $i$, $\phi_i(x,d_i)$ \th-divides over $A$. Let $c_1, \dots , c_n$ be the conjugates of $c$ over $B$ and for $0\leq j\leq n$ let $\sigma_j$ be some $B$-automorphism such that $\sigma_j(c)=c_j$. This means that $ tp(a/B)\models \vee_{j=1}^n \psi(x,c_j) $ and for any such $j$ $ \psi_j(x, c_j) \Rightarrow \vee_{i=1}^m \phi_i(x,\sigma_j(d_i)). $ Combining this two, we get \[ tp(a/B)\models \bigvee _{j=1}^n \bigvee_{i=1}^m \phi_i(x,\sigma_j(d_i)) \] and by definition $tp(a/B)$ \th-forks over $A$.
\end{proof}

We can also prove without any assumptions on the theory a version
of partial left transitivity:

\begin{lemma}
\emph{Partial Left Transitivity:} If $a\thind _A c$ and $b\thind
_{Aa} c$ then $ab\thind _A c$. \label{transitivity}
\end{lemma}

\begin{proof}
We begin with a claim.
\begin{claim}
It is enough to show that $tp(ab/Ac)$ does not \th-divide over
$A$.
\end{claim}

\begin{proof}
Suppose $a\thind_A c$, $b\thind_{Aa} c$ and $ab \nthind_A c$.
Choose $\phi (x,y,c)\in tp(ab/Ac)$ which implies $ \bigvee_{i<n}
\psi _i (x,y,c_i )$, where each $\psi_i$ \th-divides over $A$. Put
$\bar c:=(c_1:i<n)$. By extension we can find some $a'\models
tp(a/Ac)$ such that $a'\thind _A c\bar c$. Let $b_2$ be the image
of $b$ under an automorphism fixing $Ac$ which sends $a$ to $a'$,
so that $tp(ab/Ac)=tp(a'b_2/Ac)$ and $b_2\thind _{Aa'} c$. By
extension again we can also find some $b'\models tp(b_2/Aa'c)$
such that $b'\thind _{Aa'} c\bar c$. So we have tuple $a'b'$
satisfying the hypothesis of the theorem. However, $\phi
(x,y,c)\in tp(a'b'/Ac)$ and $\phi (x,y,c)\Rightarrow \bigvee \psi
_i (x,y,c_i )$, so $\psi _i (x,y,c_i )\in tp(a'b'/Ac\bar c)$ for
some $i<n$, and $tp(a'b'/Ac\bar c)$ \th-divides over $A$.
\end{proof}

Let us suppose then that $tp(ab/Ac)\models \phi (x,y,c) $ which
\th-divides over $A$ strongly dividing, say, over $Ad$. By base
extension, we can find some $d'\models tp(d/Ac)$ such that
$a\thind _{Ad'} c$, so we will assume $a\thind _{Ad} c$. By
definition of strong dividing we know that $\{ \phi(x,y,z) \}
_{z\models tp(c/Ad)}$ is $k$-inconsistent and $tp(c/Ad)$ is non
algebraic; by \ref{introduction}$(4)$ so is $tp(c/Ada)$. This
means that \[ \{ \phi (a,y,z)\} _{z\models tp(c/Ada)} \] is still
$k$-inconsistent and $tp(c/Ada)$ is non-algebraic. By definition
$tp(b/Aac)$ strongly divides over $Aad$ which means it \th-divides
over $Aa$, a contradiction.
\end{proof}

As in simple theories, this new notion of bifurcation does have
some relation with independent sequences which will be stated in
the following lemma. However, unlike in simple theories this
relation does not seem to be fundamental for the development of
the theory, perhaps because the main function of indiscernible
sequences in simple theories is to provide some uniformity when
witnessing division; such uniformity is provided by the
definability of the parameters in the the definition of
\th-dividing.

\begin{lemma}
Let a,b be elements and $A$ a set. Let $p\left( x,b\right)
=tp\left( a/Ab\right) .$ Then the following conditions are
equivalent:

\begin{enumerate}
\item  $tp(a/Ab)$ does not \th-divide over $A$

\item  For any $ B\supseteq A$  such that $b$ is not algebraic
over $B,$ there is some tuple   $a' \models tp\left( a/Ab\right) $
and some infinite $Ba'$-indiscernible sequence $I$ containing $b$.
\end{enumerate}
\end{lemma}

\begin{proof}
$\left( \Leftarrow \right) $ We will proceed by contradiction.
Assume $tp(a/Ab)$ \th-divides over $A$. By definition we have a
$B$ with $b\notin acl(B)$ and a $\delta $, such that $\models
\delta \left( a,b\right) $ and $\{\delta \left( x,b'\right)
\}_{b'\models tp\left( b/B\right) }$ is   $k$-inconsistent. Let
$I$ be any $B$-indiscernible sequence. Since the underlying set of
$I$ is a subset of $\{b'\mid b'\models tp\left( b/B\right) \}$, we
have that $\{ \delta \left( x,b'\right) \} _{b'\in I}$ is
$k$-inconsistent. Let $a'\models tp\left( a/Ab\right)$, $a'\models
\delta (x,b)$. By $k$-inconsistency we have $\nvDash \delta \left(
a',b' \right) $ for all but finite $b'\in I$ which implies that
$I$ is not $Ba'$ indiscernible.

$\left( \Rightarrow \right) $ Let $A,b,a$ be such that $tp(a/Ab)$
does not \th-divide over $A$. If $b$ is algebraic over $A$,
condition 2 in the lemma holds immediately so there is nothing to
prove. We will assume then that $b$ is not algebraic over $A$.

Let $B$ be any set containing $A$ with $b\notin acl(B)$ and let
$q\left( y\right) =tp\left( b/B\right) .$ We know that for any
$\delta \left( x,b\right) \in tp\left( a/Ab\right) $ the set
$\{\delta \left( x,b'\right) \}_{b'\models tp\left( b/B\right) }$
is not $k$-inconsistent for any $k$. By compactness
$$\bigcup_{i\in \omega} q(y_i)\cup \{ p(x,y_i) : i\in \omega\}$$ is
consistent, and realized by $a',I$, say.  For any finite set
$\Delta$ of formulas we can find by Ramsey's theorem some infinite
$I_\Delta\subseteq I$ such that $I_\Delta$ is
$\Delta$-indiscernible over $Ba'$. By compactness we can find a
$Ba'$-indiscernible sequence $J$ such that $a'b'\models p(x,y)\cup
q(y)$ for any $b'\in J$. If $\sigma$ is a $B$-automorphism with
$\sigma(b')=b$ then $\sigma(a'), \sigma(J)$ will do.
\end{proof}

\subsection{Existence}

Given some independence relation, we say that such a relation
satisfies \emph{existence} if for any tuple $a$ and any set $A$,
$a$ is independent with $A$ over $A$. This is a very useful
property for simple theories. Until this point we have studied the
behavior of \th-forking in the most general context. In this
section we will only consider theories for which \th-forking
satisfies existence: given a tuple $a$ and a set $A$, $tp(a/A)$
does not \th-fork over $A$.

As we mentioned before, unlike in simple theories there is a
difference between \th-forking and \th-dividing. As we shall see
in most of the proofs, we can usually work around this problem
since extension provides a way to reduce most of the proofs down
to \th-dividing. It would be nice however to have some idea of the
relation between the parameters that we need for \th-forking and
those that are used in the corresponding \th-dividing formulas, so
we can tell how much do we have to extend a \th-forking type
before achieving \th-division. The next lemma gives a partial
answer to this question.

\begin{lemma}
If $\delta (x,a) $ is consistent and \th -forks over $A$, as
witnessed by a disjunction $\bigvee_{i=1}^n \psi (x,a_i )$ implied
by $\delta(x,a)$, such that $\psi_i (x,a_i) $ strongly divides
over $Ac_i$, then  $a_i$ is algebraic over $Aa\bar{c}$ where $\bar
c =(c_i:i<n)$. Even more, for at least one $i$, $a_i$ is algebraic
over $Aac_i$. (We are assuming that there are no ``extra''
$\psi_i$'s: i.e. that $\delta(x,a)\nRightarrow \bigvee_{j\in I}
\psi (x,a_{i_j} )$ for any $I\subsetneq \{ 1,2,\dots, n\}$.)
\end{lemma}

\begin{proof}
Let $b$ be any element such that $b \models \delta (x,a) $. By
existence and definition of \th-forking we can extend $tp(b/Aa)$
to $tp(b/Aa\bar{a_i}\bar{c})$ so that $b\thind _{Aa}
\bar{a_i}\bar{c}$. In particular, by partial transitivity $b\thind
_{Aac_i} a_i$

We know that $b\models \psi_i (x,a_i)$ for some $i$; for any $a_i'
\models tp(a_i/Abac_i)$ we have $\mathcal{C} \models \psi(b,a_i')
$ and $a_i' \models tp(a_i/Ac_i ) $. By the definition of strong
dividing there cannot be infinitely many such $a_i'$'s ($b$ would
witness the consistency). Thus $tp(a_i/Abc_i)$ must be algebraic
and therefore so is $tp(a_i/Abac_i)$. But we know that $b \thind
_{Aac_i} a_i$, so by lemma \ref{introduction}$(4)$ this means that
$tp(a_i/Aac_i) $ is algebraic.

To finish the proof, we just have to be careful when extending
$tp(b/Aa)$ to $Aaa_i$;  let $\delta(x,a,a_i)$ be $\delta
(x,a)\wedge \neg \psi_i(x,a_i)$ so $\delta(x,a,a_i)$ implies the
disjunction $\bigvee_{j\neq i} \psi_j(x,a_j)$. Since
$k$-inconsistency is preserved, either there is some $a_j$
algebraic over $Aa_i$ or $\delta(x,a,a_i)$ \th-forks over $Aa_i$
in which case we can repeat the process and get some $a_j$
algebraic over $Aa_ic_j$. Either way we get $a_j$ algebraic over
$Aa_ic_j$ for some $j$. However, $a_i$ is algebraic over $Ac_i$ so
$a_j$ is algebraic over $Ac_ic_j$.
\end{proof}

\begin{theorem}
\label{leftextension} Let $p(x,b)$ be a type over $Ab$ which is
non-\th-forking over $A$ and let $\theta (x,y)$ be a formula such
that $p(x,b) \models \exists y \theta (x,y)$. Then $p(x,b)\cup \{
\theta (x,y) \} $ does not \th-fork over $A$.
\end{theorem}

\begin{proof}

Suppose that it does \th-fork over $A$ so that $p(x,b)\cup \{
\theta (x,y) \} \models \bigvee_{i=1}^n \psi _i (x,y,a_i)$ where
$\psi _i (x,y,a_i)$ \th-divides over $A$ by strongly dividing over
$Ac_i$.

By hypothesis $p(x,b)$ does not \th-fork over $A$, so using
extension we can choose some $a\models p(x,b)$ such that $a\thind
_A b\bar{a_i}\bar{c_i}$ which implies by partial transitivity that
$a\thind _{Ac_i} a_i$. On the other hand, using existence and
extension we can choose some $c\models \theta (a,y)$ such that
$c\thind _{Aa} b\bar{a_i}\bar{c_i}$. Using partial transitivity
again we get $c\thind _{Aac_i} a_i$.

We know that $(a,c)\models p(x,b)\cup \theta \{ (x,y)\} $ so
$(a,c)\models \psi_j(x,y,a_j)$ for some $j$. By definition, $a_j$
is algebraic over $Aacc_j$. However, using lemma
\ref{introduction} with $c\thind _{Aac_j} a_j$ we get that $a_j$
is algebraic over $Aac_j$. Using \ref{introduction} one more time
will give us $a_j$ algebraic over $Ac_j$, contradicting the
definition of strong dividing.

\end{proof}

\begin{definition}
A notion of independence $\ind ^0$ has the \emph{strong extension
property} if whenever $a \ind^0_A b$ then for any $c$ there is
there is $c'\models tp(c/Aa) $ with $ac'\ind^0 _A b$, or
equivalently, for any $c$ there is $b'\models tp(b/Aa)$ with
$ac\ind^0_Ab'$.
\end{definition}

\begin{corollary}
\th-independence has the strong extension property.
\end{corollary}

\begin{proof}

It is enough to prove that if $p(x,b) $ is a non forking extension
of $p\upharpoonright A$ and $q (x,y)$ is a consistent type over
$A$ containing $p\upharpoonright A$ then $p(x,b)\cup q(x,y)$ does
not \th-fork over $A$. This follows from theorem
\ref{leftextension}, as consistency of $(p\restriction A) \cup q$
implies consistency of $p\cup q$.
\end{proof}

\begin{definition}
A sequence of elements is a \emph{\th-Morley sequence over $A$} if
it is indiscernible and \th-independent over $A$.
\end{definition}

\begin{claim}
Let $p(x)$ be a complete type over $B\supset A$ which does not
\th-fork over $A$. Then there is a \th-Morley sequence over $A$
with all of its elements realizing $p(x)$.
\end{claim}

The claim is true for any independent notion that has extension.
The proof is exactly the same as the one given in the simple
theoretic context (see \cite{Wagner}). We will prove later that in
the theories that will actually interest us, \th-Morley sequences
will provide an alternative definition of \th-forking.

\section{\th-Rank and Rosy Theories} \label{ranksection}
\subsection{Definition of \th-rank}

We will now define a notion of rank that will code \th-forking.

\begin{definition}
Given a
formula $\phi ,$ a set $%
\Delta $ of formulas in the variables $x;y$, a set of formulae
$\Pi$ in the variables $y;z$ (with $z$ possibly of infinite
length) and a number $k,$ we define \th $ \left( \phi ,\Delta, \Pi
,k\right) $ inductively as follows:
\begin{enumerate}
\item  \th $\left( \phi ,\Delta , \Pi, k\right) \geq 0$ if $\phi $
is consistent.

\item  For $\lambda $ limit ordinal, \th $ \left( \phi ,\Delta,
\Pi ,k\right) \geq \lambda $ if and only if \th $ \left( \phi
,\Delta , \Pi, k\right) \geq \alpha $ for all $\alpha <\lambda $

\item  \th $ \left( \phi ,\Delta, \Pi ,k\right) \geq \alpha +1$ if
and only if there is a $ \delta \in \Delta $, some $\pi(y;z) \in
\Pi$ and parameters $c$ such that

\begin{enumerate}
\item  \th $ \left( \phi \wedge \delta \left( x,a\right) ,\Delta
,\Pi ,k\right) \geq \alpha $ for infinitely many $a\models
\pi(y;c) $

\item  $\left\{ \delta \left( x,a\right) \right\} _{a \models
\pi(y;c)}$ is $k-$inconsistent
\end{enumerate}
\end{enumerate}

\end{definition}

As usual, for a type $p$, we define
\[ \text{\th}  \left( p,\Delta ,\Pi ,
k\right) =\min \left\{ \text{ \th } \left( \phi \left( x\right)
,\Delta ,\Pi ,k\right) \mid \phi \left( x\right) \in p\right\} .\]

\begin{remark}
We can give a definition directly for types, changing all
instances of $\phi$ for some type $p$ which is the case we will
usually use. However, when dealing with types, we will use an
alternative version of condition 3(a). Let $p$ be a complete type
over $A$ and let us assume that we are witnessing the rank going
up as in condition 3, with $c$ being a tuple containing all the
parameters in $\pi$. We can always find some non-algebraic
complete type $q(y)$ over $Ac$ containing $\pi$ such that for any
$a'\models q(y)$, \th $ \left( \phi \wedge \delta \left( x,a'
\right) ,\Delta ,\Pi ,k\right) \geq \alpha $. Therefore, we can
change condition 3(a) by

3(a)'. \th$\left( p\cup \{\delta(x,a) \}, \Delta, \Pi, k\right)
\geq \alpha$ where $a\models \pi$ and $tp(a/Ac)$ is non-algebraic.

\end{remark}

\begin{remark}
\label{typedefinability} Given any type $p(x;y)$ (not necessarily
complete), formulas $\phi, \pi$ and integers $k,n$, the set

\[ \left\{ b | \text{\th} \left( p(x,b), \phi, \pi, k \right) \geq n \right\}\] is type definable.

This has nice consequences for the structure of the \th-rank. For
example, compactness implies that \th$\left( p\cup \{\delta(x,a)
\}, \Delta, \Pi, k\right)$ is finite whenever it is defined.

\end{remark}

\begin{proof}
Notice that \[ \text{\th} (p(x,b), \phi, \pi, k)\geq n \] is
witnessed by a tree where each of the nodes is a formula
$\phi(x,a)$, its height is $n$, for each $i\leq n$ there is some
$c_i$ such that the level $i$ contains all $\phi(x,a)_{a\models
\pi(x,c_i)}$, the union any $k$ of such formulas (all in the same
level) is inconsistent and the union of any branch is consistent
with $p$. All theses properties can be described by formulas, and
we get type definability.
\end{proof}

\subsection{Properties of the \th-rank}

\begin{theorem}
\label{propertiesranks} This thorn rank has the following
properties:

\begin{enumerate}
\item  \emph{Monotonicity:} If $\Delta \subseteq \Delta '$,
$p\supseteq p',$ and $\Pi \subseteq \Pi'$ then
\[
\text{\th} \left( p',\Delta ', \Pi ',k\right) \geq \text{ \th}
\left( p,\Delta ,\Pi ,k\right).
\]

\item  \emph{Transitivity:} If $p\subseteq q\subseteq r,$ then \th
$ \left( r,\Delta, \Pi ,k\right) =\text{ \th} \left( p,\Delta ,\Pi
,k\right)$ if and only if $  \text{ \th} \left( r,\Delta ,\Pi
,k\right) =\text{ \th} \left( q,\Delta ,\Pi ,k\right) $ and $
\text{ \th} \left( q,\Delta ,\Pi ,k\right) =\text{ \th} \left(
p,\Delta ,\Pi ,k\right) $.

\item  \emph{Additivity:} \th $ \left( \left( \theta \vee \psi
\right) ,\Delta , \Pi ,k\right) =\max \left\{ \text{ \th} \left(
\theta ,\Delta ,\Pi ,k\right) ,\text{ \th} \left( \psi ,\Delta
,\Pi ,k\right) \right\} $
\end{enumerate}
\end{theorem}

\begin{proof}.

\noindent \emph{1)} We will prove by induction on $\lambda $ that
if $\lambda \leq \text{ \th} \left( p,\Delta , \Pi ,k\right) $
then \th $ \left( p ',\Delta ', \Pi' , k\right) \geq \lambda .$
For $\lambda =0$, the proof is clear, as is the induction step for
the case when $ \lambda $ is a limit ordinal. Now, suppose it is
true for $\lambda=\alpha $ and let us assume that \th $\left(
p,\Delta ,\Pi ,k\right) \geq \alpha +1$. Let $\psi$ be any formula
in $p'$; by hypothesis $\psi\in p$ so \th $\left( \psi ,\Delta
,\Pi ,k\right) \geq \alpha +1$. We can therefore find a $\phi \in
\Delta $ and a formula $\pi $ in $\Pi$ such that \th $ \left( \psi
\wedge \phi \left( x,a\right) ,\Delta ,\Pi ,k\right) \geq \alpha $
for infinitely many $a\models \pi$ and $\left\{ \phi \left(
x,a\right) \right\} _{a\models \pi }$ is $k$ -inconsistent$.$ By
induction hypothesis \th $ \left( \psi \wedge \phi \left(
x,a\right) ,\Delta ', \Pi' ,k\right) \geq \alpha $. But $\psi \in
p'$, $\phi \in \Delta '$ and $\pi\in \Pi'$ so $\phi $ actually
witnesses \th $\left( \psi ,\Delta ', \Pi', k\right) \geq \alpha
+1$ for all $\psi \in p'$.

\noindent \emph{2)} By monotonicity, \th$(r, \Delta, \Pi, k)\leq
\text{\th} (q, \Delta, \Pi, k)\leq \text{\th}(r, \Delta, \Pi,
k)\geq $; transitivity follows from transitivity for equality.

\noindent \emph{3)} By monotonicity we have
\[ \text{ \th}
\left(  \left( \theta \vee \psi \right) ,\Delta ,\Pi , k\right)
\geq \max \left\{ \text{ \th} \left(  \theta ,\Delta , \Pi
,k\right) , \text{ \th} \left(  \psi ,\Delta , \Pi , k\right)
\right\} . \]

For the other direction we will prove by induction on $\alpha $
that if
\[ \text{\th} \left(  \left( \theta \vee \psi \right)
,\Delta ,\Pi ,k\right) \geq \alpha \] then either \th $ \left(
 \theta ,\Delta ,\Pi ,k\right) \geq \alpha $ or \th $
\left(  \psi ,\Delta ,\Pi ,k\right) \geq \alpha .$ Once again, the
only difficult step is the induction step. Let us assume it is
true for $\alpha $ and that \th $ \left(  (\theta \vee \psi
),\Delta , \Pi ,k\right) \geq \alpha +1.$ We can then find a
$\delta \in \Delta $ and a formula $\pi \in \Pi $ such that
\[ \text{ \th}
\left(  (\theta \vee \psi )\wedge \delta \left( x,a_{i}\right)
,\Delta ,\Pi ,k\right) \geq \alpha \] for infinitely many
$a_{i}\models \pi $ and $\left\{ \phi \left( x,a_{i}\right)
\right\} _{a_{i}\models \pi }$ is $k$-inconsistent. By induction,
for any such $a_{i}$ either \th $ \left(  \theta \wedge \delta
\left( x,a_{i}\right) ,\Delta , \Pi ,k\right) \geq \alpha $ or \th
$ \left(  \psi \wedge \delta \left( x,a_{i}\right) ,\Delta ,\Pi
,k\right) \geq \alpha$; but then for one of $\psi $ or $\theta $
(let us assume $\theta $) we have infinitely many $a_{i}'s$ such
that \th $ \left(  \theta \wedge \delta \left( x,a_{i}\right)
,\Delta ,\Pi ,k\right) \geq \alpha $ and by definition \th $
\left(  \theta ,\Delta ,\Pi ,k\right) \geq \alpha +1.$
\end{proof}

\begin{corollary}{\emph{Extension for fixed $\Delta $ and $\Pi$:}}

\noindent For any partial type $\xi $ defined over a set $A$,
finite sets of formulas $\Delta $ and $\Pi$ and any $k,$ we can
extend $\xi $ to a complete type $p$ over $A$ such that \th $
\left( p,\Delta , \Pi , k\right) =\text{\th} \left( \xi ,\Delta
,\Pi ,k\right) .$
\end{corollary}

\begin{proof} It is just the usual application of additivity, the
definition of \th $ $-rank for types and Zorn's lemma as it is
used in simple theories (see ~\cite{Kim}).
\end{proof}

One of the properties that we will prove here is that
\th-independence and \th-ranks are related in the same way that
forking and the $D$-ranks are. The next theorem is one of the
directions of the relation we will prove.

\begin{theorem}
\label{d-ranks} Let $p$ be a type over $B\supseteq A$ such that
for every $\phi ,\Pi ,k,$
\[
\text{\th} \left( p\upharpoonright A,\phi ,\Pi ,k\right) =
\text{\th} \left( p,\phi ,\Pi ,k\right)\footnote{ By $\text{\th}
\left( p,\phi ,\Pi , k\right)$ we really mean $\text{\th} \left(
p,\{\phi\} ,\Pi , k\right)$.}.
\] Then $p$ does not \th-fork over $A.$
\label{ranks}
\end{theorem}

\begin{proof}
Suppose $p$ \th-forks over $A$. So
\[p\models \bigvee _{i<n}\phi
_{i}\left( x,b_{i}\right) \] where each $\phi _{i}$ \th-divides
over $A$. By definition there is some $\theta_i(y,z)$, an element
$d_i$ and some finite $k_i$ such that $\{\phi (x,b_i')
\}_{b_i'\models \theta_i(y,d_i)}$ is $k_i$-inconsistent and
$tp(b_i/Ad_i)$ is a non-algebraic type containing $\theta(y,d_i)$.
We can always take the maximum of the $k_i$'s and assume they are
all equal.

Now, as in the proof for simple theories (see \cite{Kim}), we
construct a formula
\begin{equation*}
\psi \left( x,y_{1},\dots ,y_{n},z\right) :=\bigvee_{i<n}\left(
\phi _{i}\left( x,y_{i}\right) \wedge y_{i}=z\right)
\end{equation*}
which uniformizes the $\phi _{i\text{ }}$ in the following sense:
if $ \bar{y}=\left( y_{1},\dots ,y_{n},w\right) $ then for any
$i\leq n$,  $w=y_{i}$ implies $ \psi \left( x,\bar{y}\right) \iff
\phi _{i}\left(x, y_{i}\right) .$ On the other hand, if we define
a new formula $\theta _{i}' \left( \bar{y}, z\right)
:=\{w=y_{i}\}\cup \theta _{i}\left( y_{i}, z\right) $ we have that
for any tuple $\bar{b}=\left< b_1, b_2, \dots , b_n,l \right>$,
and for any $d$, $\bar{b}\models \theta _{i}'(y,d)$ if and only if
$ b_{i}\models \theta _{i} (y,d)$ and $l=b_i$. Thus, if we define
$c_{i}:=\left< b_{1},\dots ,b_{n},b_{i}\right> $ we get $p\models
\vee _{i<n}\psi \left( x,c_{i}\right) $ and $\{\psi \left(
x,c_{i}\right) \}_{c_{i}\models \theta _{i}'(y,d) }$ is
$k$-inconsistent. By additivity we can add at least one of the
$\phi _i$'s without changing the rank, so we can extend $p$ to
some type $q$ over $Bc_{1}\dots c_{n}$ such that \th $\left(
p,\psi ,\{\theta_{1}' ,\dots \theta_{n}'\},k\right) $=\th $\left(
q,\psi ,\{\theta _{1}' ,\dots \theta_{n}' \},k\right)$ and $q$
implies one of the $\phi _{i}\left( x,b_{i}\right) $ and thus one
of the $\psi \left( x,c_{i}\right) .$ Whichever one it is (we can
assume it is the first one without loss of generality), we know
that $\{\psi \left( x,t\right) \}_{t\models \theta _{1}'(y,b) }$
is $k$-inconsistent and $tp(\bar{b} /Ad)$ is a non-algebraic type
containing $\theta'_1(y,d)$. By definition of \th-rank
\[ \text{\th} \left( p\upharpoonright A,\psi ,\{\theta' _{1},\dots \theta' _{n}
\},k\right) \geq \text{\th} \left( q,\psi ,\{\theta' _{1} ,\dots
\theta' _{n} \},k\right) +1\] and therefore
\begin{equation*}
\text{\th }\left( p,\psi ,\{\theta' _{1},\dots \theta' _{n}
\},k\right) <\text{\th }\left( p\upharpoonright A,\psi ,\{\theta'
_{1},\dots \theta' _{n} \},k\right),
\end{equation*}
a contradiction.
\end{proof}

\begin{corollary}
If for any finite $\Delta$, $\Pi$ all the \th-ranks are defined,
then we have existence (and all the results of subsection
\ref{ranksection}).
\end{corollary}

\begin{proof}
Let $p(x)$ be a type over $A$. Then $p\restriction A=p$ and for
any finite set of formulas $\Delta$, $\Pi$ and any finite $k$,
\[ \text{\th} \left( p\upharpoonright A,\phi ,\Pi ,k\right) =
\text{\th} \left( p,\phi ,\Pi ,k\right).\] The theorem implies
that $p(x)$ does not \th-fork over $A$.
\end{proof}

\section{Rosy Theories}

As the corollary above suggests, there is a lot to say about
theories that have ordinal-valued \th-ranks even if we limit
ourselves to finite formulas and types in the definition. Doing
this has the extra advantage that, if we limit ourselves to finite
$\Delta$ and $\Pi$, definability of \th-rank and compactness give
us that a \th-rank is defined if and only if it is finite.

From now on, we will study the class of theories such that
\th$(p,\Delta, \Pi,k)$ is finite for any type $p(x)$ in the
language of the theory, any finite sets of formulas $\Delta$ and
$\Pi$ and any finite number $k$. We will call any such theory
``rosy''\footnote{Rosy theories include simple and o-minimal
theories but, as mentioned in the introduction, there are other
rosy theories which do not fall into either of this categories.}.
We will prove that in any rosy theory \th-forking has a lot of the
geometric properties we want for an independence notion. All
theories in this section are assumed to be rosy.

\subsection{Geometry of \th-independence in Rosy theories}

\begin{theorem} \emph{Symmetry:} For any two elements $a,b$ and any set $A$, $a\thind _A b$ if and only if $b\thind _A a$.
\end{theorem}

\begin{proof}

Let us suppose $a\thind_A b$ and let us assume that $\models
\delta (b,a)$ where $\delta(x,a)$ \th-forks over $A$. Let $\delta
(x,a) \Rightarrow \bigvee \psi_i (x,c_i) $ where each $ \psi_i
(x,c_i) $ strongly divides over $Ad_i$. By compactness we can
prove the $k$-inconsistency using some formula $\pi_i$ with
parameters in $Ad_i$.

We define a sequence $\left\{ a^j, \bar{c^j}, \bar{d^j}
\right\}_{j\in \mathbb{N}}$ in the following way. Let $a^0=a$,
$\bar{c^0}=\bar{c_i}$ and $ \bar{d^0}=\bar{d_i}$. Assuming we have
defined the sequence up to $j=n$, let $a^{n+1}\models tp(a/Ab)$ be
such that

\[ a^{n+1}\thind _A
ba^1\bar{c_i}^1\bar{d_i}^1 \dots a^n\bar{c_i}^n\bar{d_i}^n
\] (we know such $a^{n+1}$ exists by extension)
and let $c_i^{n+1}$ and $d_i^{n+1}$ be the images of $c_i$ and
$d_i$ under an automorphism that sends $a$ to $a^{n+1}$.

Once we have such a sequence, by right partial transitivity and
monotonicity (\ref{introduction}) we have that for all $n>m$,
\[ a^n\thind _{Ad_i^ma^{m+1}\dots a^{n-1}} c_i^m.\]
Using induction and left partial transitivity (\ref{transitivity})
$n-m$ times we get
\[
a^{m+1}a^{m+2}\dots a^n \thind _{Ad_i^m} c_i^m
\] so in particular it cannot strongly divide. But $tp(c_i^m/Ad_i^m)$
is non-algebraic so
\[ tp(c_i^m/d_i^ma^{m+1}a^{m+2}\dots a^n)\] is
non-algebraic for any $n>m$; by local character of \th-forking and
compactness
\[ tp(c_i^m/d_i^ma^{m+1}a^{m+2}\dots )\]
would also be non-algebraic. By definition, $\psi _i (x,c_i^n)$
strongly divides over $Ad_i^na^{n+1} \dots$ and $\psi _i
(x,c_i^n)$ \th-divides over $Aa^{n+1}a^{n+2}\dots $.

Now, since $a^n\bar{c_i}^n \models tp(a\bar{c_i}/A)$, $a^n\models
tp(a/Ab)$ we have that for all $n$
\[\delta (x,a^n) \Rightarrow \bigvee \psi _i (x,c_i^n)\] and $b \models
\delta (x,a^n)$. This implies $\delta(x,a^n)\in
tp(b/Aa^na^{n+1}a^{n+2}\dots )$ \th-forks over

\noindent $Aa^{n+1}a^{n+2}\dots $. By the proof of the previous
theorem we know that this is witnessed by
\begin{equation*}
\text{\th }\left( tp\left( b/Aa^{n}a^{n+1}a^{n+2}\dots \right)
,\psi ,\theta ,k\right) <\text{ \th }\left( tp\left(
b/Aa^{n+1}a^{n+2}\dots \right) ,\psi ,\theta ,k\right)
\end{equation*}
where $\theta $ is the one mentioned at the end of the proof of
theorem \ref{ranks} and $\psi $ depends only on the $\psi _{i}$
(so $\theta$ and $\psi_i$ are the same for all $n$). This means
that \th$\left( tp(b/Aa^1a^2a^3\dots ) ,\psi ,\theta ,k\right)$ is
infinite contradicting rosiness.
\end{proof}

If we combine this result with those above, we get an analogue for
most of the geometric properties of forking independence for
simple theories (see ~\cite{Wagner} 2.3.13). We will use the same
names.

\begin{corollary}
\th-independence defines an independence relation --in the sense
described in \cite{Kim/Pillay}-- in any rosy theory. More
precisely, if we are working inside a theory $T$ of finite
\th-ranks, then the following properties hold for types in models
of $T$.
\begin{enumerate}

\label{allproperties} \item \emph{Existence:} If $p\in S(A)$, then
$p$ does not \th-fork over $A$. \item \emph{Extension:} Every
partial type over $B\supset A$ which does not \th-fork over $A$
can be extended to a complete type $p(x)$ over $B$ which does not
\th-fork over $A$. \item \emph{Reflexivity:} $B\thind _A B$ if and
only if $B\subseteq acl(A)$. \item \emph{Monotonicity:} If $p$ and
$q$ are types with $p\supseteq q$ and $p$ does not \th-fork over
$A$, then $q$ does not \th-fork over $A$. \item \emph{Finite
Character:} $C\thind _A B $ if and only if $c\thind _A B$ for any
finite $c\subset C$. \item \emph{Symmetry:} $C\thind _A B$ if and
only if $B\thind_A C$. \item \emph{Transitivity:} If $A\subset
B\subset C$ then $b\thind _A C$ if and only if $b\thind _A B$ and
$b\thind _B C$. \item If $a\thind_{A} b$ and for some
formula $%
\delta(x,y)$  $\delta \left( x,a\right) $ \th -forks over $A$,
then $\delta \left( x,a\right) $ \th -forks over $Ab.$ \item Let
$A\subset B$. If $a\thind _A B$ then $a\thind _A acl(B)$.
\end{enumerate}
\label{properties}
\end{corollary}

\begin{proof}
Properties
\emph{1},\emph{2},\emph{3},\emph{4},\emph{6},\emph{7},\emph{8} and
\emph{10} have already been proven or follow immediately as a
corollary of symmetry. 5. is clear from the definitions.

To prove 9. we will first show that if $\delta (x,a)$ \th-divides
over $A$, then it \th-divides over $Ab$. Take $B\supseteq A$ such
that $\left\{ \delta \left( x,a'\right) \right\} _{a'\models
tp\left( a/B\right) }$ is $k$-inconsistent and $tp(a/B)$ is
non-algebraic. By base extension, we can find some $B'\models
tp(B/Aa)$ such that $a\thind _{B'} b$; note that $\delta (x,a)$
strongly divides over $B'$ so $tp(a/B')$ is non-algebraic.  By
\ref{introduction} $tp\left( a/B'b\right) $ is not algebraic and $
\left\{ \delta \left( x,a'\right) \right\} _{a'\models tp\left(
a/B'b\right) }$ is still $k$-inconsistent. By definition $\delta
(x,a)$ \th-divides over $Ab$.

Suppose now that $\delta \left( x,a\right) $ \th -forks over $A$
and $\delta \left( x,a\right) \Rightarrow \bigvee \phi _{i}\left(
x,a_i\right) $ where each formula $\phi _i(x,a_i)$ \th-divides
over $A$. By extension we can assume $a_ia\thind _A b$ which
implies $a_i\thind _A b$. But we just proved that $\phi _i
(x,a_i)$ \th-divides over $Ab$ and thus $\delta (x,a)$ \th-forks
over $Ab$.

\end{proof}

As we mentioned before, in any rosy theory \th-forking is
witnessed by the \th-ranks in the same way that forking was
witnessed by the $D$-ranks. We already proved one of the
implications in \ref{d-ranks}; the following theorem will prove
the other one.

\begin{theorem}
Let $p$ be a type over $B\supseteq A$ such that $p$ does not \th
-$fork$ over $ A$. Then for all $\phi ,\theta ,k$ \th $\left(
p\upharpoonright A,\phi ,\theta ,k\right) $=\th $\left( p,\phi
,\theta ,k\right) .$ \label{ranks2}
\end{theorem}

\begin{proof}
Since whenever $\Delta$ and $\Pi$ are finite \th$(p,\Delta, \Pi,
k)$ is finite for any type $p$ and any $k$, it is enough to prove
that for any natural number $m$, \th $\left( p\upharpoonright
A,\phi ,\theta ,k\right) \geq m$ implies \th $\left( p,\phi
,\theta ,k\right) \geq m.$ The case $m=0$ is clear by the
consistency of $p.$ Suppose then that we have it for $m=n$ and
that \th $\left( p\upharpoonright A,\phi ,\theta ,k\right) \geq
n+1.$ By definition we can find some $b$ and some tuple $c$ such
that \th $\left( p\upharpoonright A\cup \{\phi \left( x,b\right)
\},\phi ,\theta ,k\right) \geq n$, $b\models \theta(y,c)$, $\{\phi
\left( x,b'\right) \}_{b'\models \theta(y,c) }$ is $k$
-inconsistent and $tp\left( b/Ac\right) $ is non-algebraic.

By \ref{propertiesranks} we can extend $p\upharpoonright A\cup
\{\phi \left( x,b\right)\}$ to a type $ q\left( x,b\right) $ over
$Ab$ such that
\[\text{\th} \left(
q\left( x,b\right) ,\phi ,\theta ,k\right) \geq n.
\]

Now, if $a\models p$ and $a'\models q$, we have an automorphism
that fixes $A$ and sends $a$ to $a'$, $p$ to some type $p'$ over
$B'$ which does not \th-fork over $A$ and $c$ to some $c'$ which
still witnesses the \th-division of $\phi(x,b)$. Since $\phi$ and
$\theta$ have no parameters they are fixed under automorphisms so
it is enough to show that \th $(p',\phi, \theta ,k)\geq n +1$. We
may therefore assume without loss of generality that $p=p',$
$a=a'$ and thus $a$ satisfies both $q$ and $p.$ By extension we
can extend $ tp(b/Aa)$ to a type over $Ba$ which does not \th-fork
over $Aa$; let $b'$ be some realization of such extension. Let
$\sigma$ be some automorphism that sends $b$ to $b'$, fixes fixing
$Aa$ and sends $c$ to $c'$.

We may assume $b'=b$ and $b\thind _{Aa}B$. By hypothesis,
$a\mathop{\mathpalette\Ind{}}_{A}^{\text{\th }}B$ and by
transitivity, $tp\left( ba/B\right) $ does not \th -fork over $A$.
This implies that $tp\left( a/Bb\right) $ does not \th -fork over
$Ab$ (by transitivity) and $tp\left( B/baA\right) $ does not \th
-fork over $A$ by symmetry of \th-forking. From this last result,
using transitivity again, we get that $B\thind_{A} b$ and $B
\thind_{Ab} a$. By symmetry $ b\thind_{A} B$ and $tp\left(
a/Bb\right) :=q_{B}\left( x,b\right) $ does not \th-fork over
$Ab.$ But $tp\left( a/Ab\right)$ is $q(x,b)$ so by induction
hypothesis
\begin{equation*}
\text{\th }\left( tp\left( a/Bb\right) ,\phi ,\theta ,k\right)
=\text{\th } \left( q\left( x,b\right) ,\phi ,\theta ,k\right)
\geq n.
\end{equation*}

Since $b \thind_{A}B$, we know we can find some ${c'}\models
tp({c}/Ab)$ such that $b\thind_{A{c'}}B$, $\models \theta
(b,{c'})$ and $\{ \phi (x,y)\} _{\models \theta (y,{c'})}$ is
$k$-inconsistent. These are the only things we need ${c}$ for so
we will assume $b\thind _{A{c}} B$. By hypothesis $tp\left(
b/A{c}\right) $ is non-algebraic so by lemma \ref{introduction}
neither is $tp\left( b/B{c}\right) .$

Now, for any $b'\models tp\left( b/B{c}\right) $ we have \th
$\left( q_{B}\left( x,b'\right) ,\phi ,\theta ,k\right) \geq n$,
$\phi \left( x,b'\right) \in q_{B}\left( x,b'\right) $ and $\{\phi
\left( x,b'\right) \}_{b'\models \theta \left( y,{c}\right) }$ is
$k$-inconsistent. Therefore,
\begin{equation*}
\text{\th }\left( q_{B}\left( x,b'\right) ,\phi ,\theta ,k\right)
\geq n
\end{equation*} which by monotonicity of the \th-ranks implies \begin{equation*} \text{\th }\left( p\left( x\right)
\cup \{\phi \left( x,b'\right) \},\phi ,\theta ,k\right) \geq n.
\end{equation*}
By monotonicity and definition of the \th-ranks,
\begin{equation*}
\text{\th }\left( p,\phi ,\theta ,k\right) \geq n+1
\end{equation*}
\end{proof}

\subsubsection{\th-forking and \th-Morley sequences}

Throughout the development of stability and simplicity theory
Morley sequences have been used for characterizing forking in both
stability and simplicity theory. In rosy theories we can partially
recover such characterization \footnote{Partially in the sense
that unlike the simple case, it is not true that if we take any
\th-Morley sequence the conjunction of the resulting formulas
would necessarily be consistent --o-minimal theories can witness
the failure of this remark--}.

\begin{remark}
Once we have transitivity and symmetry, any \th-Morley sequence
$<a_i>_{i\in I}$ has the property that for any $i\in I$,
$tp(a_i/A\cup \{a_j\}_{ j\in I, j\neq i})$ does not \th-fork over
$A$.
\end{remark}

\begin{proof}
The proof given for Morley sequences in simple theories is a
straightforward application of symmetry, transitivity and local
character.
\end{proof}

\begin{theorem}
If there is a \th-Morley sequence over $A$ with $a=a_{0}$ such
that $\bigwedge_i \delta \left( x,a_{i}\right) $ is consistent,
then $\delta \left( x,a\right) $ does not \th -fork over $A.$
\end{theorem}

\begin{proof}
Let $\left< a_i \right>_{i\in \mathbb{N}}$ be a \th-Morley
sequence over $A$ with $a_0=a$ and let $\bigwedge_i \delta \left(
x,a_{i}\right) $ be consistent, realized by an element so let $b$.
Suppose that $tp\left( b/Aa_{0}\right) $ \th -forks over A. By the
previous theorem we know that this is witnessed by
\begin{equation*}
\text{\th }\left( tp\left( b/Aa_{0}\right) ,\psi ,\theta ,k\right)
<\text{ \th }\left( tp\left( b/A\right) ,\psi ,\theta ,k\right)
\end{equation*}
where $\psi $ and $\theta$ are formulas which only depend on the
formulas that witness the \th-forking of $\delta(x,a)$\footnote{By
this we mean that $\psi$ and $\theta$ depend only on the formulas
which $\delta(x,a)$ implies and \th-divide over $A$ and on the
formulas needed to prove the $k$-inconsistency of such formulas.}
(see the proof of theorem \ref{ranks}).

Using an automorphism we get that $ tp\left( b/Aa_{n}\right) $ \th
-forks over $A$. However,
\[a_{n+1}\thind_{A}a_{0},\dots
,a_{n}\] so by corollary \ref{properties} we know that $\delta
\left( x,a_{n+1}\right) $ actually \th-forks over
$Aa_{0}a_{1}\dots a_{n}$ and it implies a similar disjunction of
formulas which \th-divide. By the proof of theorem \ref{ranks},
not only do we know that $tp\left( b/Aa_{0}a_{1}\dots
a_{n+1}\right) $ is a \th-forking extension of $tp\left(
b/Aa_{0}a_{1}\dots a_{n}\right) $ but we also know that for all
$n,$
\begin{equation*}
\text{\th }\left( tp\left( b/Aa_{0}a_{1}\dots a_{n+1}\right) ,\psi
,\theta ,k\right) <\text{\th }\left( tp\left( b/Aa_{0}a_{1}\dots
a_{n}\right) ,\psi ,\theta ,k\right)
\end{equation*}
which is impossible since \th $\left( tp\left( b/Aa_{0}a_{1}\dots
a_{n+1}\dots \right) ,\psi ,\theta ,k\right) $ is finite.
\end{proof}

\begin{theorem}
Suppose $\delta \left( x,a\right) $ does not \th -fork over $A.$
Then there is a \th-Morley sequence $I$ over $A$, $a_0\in I$ such
that $\bigwedge_{a_{i}\in I}\delta \left( x,a_{i}\right) $ is
consistent.
\end{theorem}

\begin{proof}
By extension, if $\delta \left( x,a\right) $ does not \th-fork
over $A,$ we can find some $b$ realizing $\delta(x,a)$ such that
$tp\left( b/Aa\right) $ does not \th -fork over $A$. By symmetry,
$tp\left( a/Ab\right) $ does not \th-fork over $A$.  We can then
construct a \th-Morley sequence $I=<a_i>$ over $A$ where $a_0=a$
and for all $i$, $a_i\models tp(a/Ab)$. Such $b$ witnesses the
consistency of
\[
\bigwedge_{a_{i}\in I}\delta \left( x,a_{i}\right). \] \end{proof}

\subsubsection{Superrosy theories and The $\uth $-rank}

As with simple theories in some cases we can define a global rank,
the $ \uth $-rank, which will share many of the properties that
the $\text{U}$-rank has. It will also help us analyze (in chapter
3) the relation between \th-forking and usual forking in
supersimple theories.

\begin{definition}
We define the $\uth$-rank inductively as follows. Let $p(x)$ be a
type over some set $A$. Then,
\begin{enumerate}
\item $\uth(p(x))\geq 0$ if $p(x)$ is consistent. \item For any
ordinal $\alpha$, $\uth(p(x))\geq \alpha+1$ if there is some tuple
$a$ and some type $q(x,a)$ over $Aa$ such that $q(x,a)\supset
p(x)$, $\uth \left( q(x,a) \right) \geq \alpha$ and $q(x,a)$
\th-forks over $A$. \item For any $\lambda$ limit ordinal,
$p(x)\geq \lambda$ if and only if $p(x)\geq \sigma $ for all
$\sigma < \lambda$.
\end{enumerate}

\end{definition}

\begin{remark}
Using extension to get from \th-forking to \th-dividing, one can
easily verify that we can replace condition (2) by the following
statement:
\begin{list}{(2')}{}
\item For any ordinal $\alpha$, $\uth(p(x))\geq \alpha+1$ if there
is some $p(x,a)\supset p(x)$ and some $c$ such that

\begin{itemize}
\item $\uth(p(x,a))\geq \alpha$. \item $tp(a/Ac)$ is
non-algebraic. \item $ \{ p(x,a') \}_{a'\models tp(a/Ac)}$ is
$k$-inconsistent for some $k$.
\end{itemize}
\end{list}
\label{remark}
\end{remark}

\begin{definition}
A theory $T$ is \emph{superrosy} if given any set $A$ in a model
of $T$ and a type $p(x)\in S(A)$, $\uth (p(x))$ is defined.

Equivalently, a theory is superrosy if there are no infinite
forking chains of types.
\end{definition}

In the case of supersimple theories (theories where the SU-rank of
any type is defined), one of the most used properties is that this
rank satisfies Lascar's inequalities. With the properties we have
proved for \th-forking, the proof of the Lascar's inequalities for
simple theories (see ~\cite{Wagner}) works in our context.

\begin{theorem}{\emph{Lascar's Inequalities:}}
Whenever the $\uth$-rank is defined, it satisfies the following
inequality:
\[ \uth(tp(a/bA))+\uth(tp(b/A)) \leq \uth (tp(ab/A)) \leq \uth
(tp(a/bA))\oplus \uth(tp(b/A)) \]
\end{theorem}
\begin{proof}
See ~\cite{Wagner} theorem 5.1.6.
\end{proof}

Finally, we want to give an alternative definition of $\uth$-rank
which is closer to the work in hyperimaginaries and because of
this it will prove very useful when comparing our $\uth$ and the
usual $\text{SU}$-rank in supersimple theories.

\begin{definition}
We define the $\uth_*$-rank on types inductively as follows. Let
$p(x)$ be a type over some set $A$. Then,
\begin{enumerate}
\item $\uth _*(p(x))\geq 0$ if $p(x)$ is consistent. \item For any
ordinal $\alpha$, $\uth _*(p(x))\geq \alpha+1$ if there is some
$p(x,a)$  and some set $B\supset A$ such that
\begin{itemize}
\item $p(x,a)\supset p(x)$. \item $\uth _*(p(x,a'))\geq \alpha $
for any $a'\models tp(a/B)$
 and $tp(a/B)$ non-algebraic.
\item There is a $k\in \mathbb{N}$ such that for any distinct
$a_1, a_2\dots , a_k \models tp(a/B)$, $\bigcup_{i=1}^k p(x,a_i)$
is a \th-forking extension of $p(x,a_1)$\footnote{If
$\bigcup_{i=1}^k p(x,a_i)$ is inconsistent it is considered to be
a \th-forking extension of $p(x,a_1)$}.
\end{itemize}
\item For any $\lambda$ limit ordinal, $p(x)\geq \lambda$ if
$p(x)\geq \sigma $ for all $\sigma < \lambda$
\end{enumerate}
\end{definition}

The main part of proving the two definitions are equivalent is the
following lemma, which in itself is a weak amalgamation for
\th-independent extensions of a given type. It is far away,
however, from the Independence Theorem that characterizes simple
theories. As mentioned in the introduction, in chapter 4 we will
talk about which of these amalgamation theorems can be true in a
general rosy theory.

\begin{lemma}
Let $A,B$ be supersets of some set $C$ such that $A\thind_C B$.
Let $p_A$ and $p_B$ be non-\th-forking extensions to $A$ and $B$
of the same type $p$ over $C$. Then there is some $B'\models
tp(B/C)$, $B'\thind _C A$ such that $p_{B'}\cup p_A$ does not
\th-fork over $C$. \label{amalgamation}
\end{lemma}

\begin{proof}
Let $a,b$ be two elements such that $a\models p_A$ and $b\models
p_B$. By hypothesis, $a\thind _C A$ and $b\thind _C B$. Since $a$
and $b$ satisfy the same type over $C$, there is an automorphism
$\sigma$ fixing $C$ such that $\sigma(b)=a$; let $B'':=\sigma(B)$
so $a\thind _C B''$. Using symmetry and extension, we can find
some $B'\models tp(B''/Ca)$ such that $B'\thind _C Aa$ (so in
particular $a\models p_{B'}$ and $B'\thind _C A$). Using
transitivity and symmetry once more we get $a\thind _AC B'$ and by
transitivity $a\thind _C AB'$.

By construction, $B'\models tp(B/C)$ and $B'\thind _C A$. Finally,
$a\models p_A\cup p_{B'}$ and $a\thind _C AB'$.
\end{proof}

\begin{claim}
For any complete type $p$,
\[
\uth (p)=\uth_* (p). \]
\end{claim}

\begin{proof}
We will prove that for any complete type $p$, $\uth(p)\geq \uth_*
(p)$ and $\uth(p) \leq \uth_* (p)$.

To prove that $\uth(p) \geq \uth_*(p)$ we will prove that for any
ordinal $\alpha$, $\uth (p)\geq \alpha$ implies $\uth_*(p)\geq
\alpha$. For $\alpha=0$ and or when $\alpha$ is a limit ordinal,
the induction follows from the definitions; the successor case
follows from remark \ref{remark}.

As for the other direction, we need to prove that

\begin{eqnarray}
\text{if } \uth_*(p)\geq \alpha \text{ then } \uth(p)\geq \alpha.
\end{eqnarray}

We will again do an induction on $\alpha$. For $\alpha=0$ or
$\alpha$ a limit ordinal the induction is immediate. For the
successor case, let us assume that (2.1) is true for all $\alpha
\leq \sigma$ and that $\uth_*(p)\geq \sigma +1$. By definition we
can find some $a$ and some $B\supset A$ such that $p(x,a)$ is an
extension of $p(x)$, $\uth_*(p(x,a)) \geq \sigma$ and there is
some $k$ such that for any $a_1, a_2, \dots a_k \models tp(a/B)$,
$\uth _*(\cup_{i=1}^k p(x,a_i)) < \uth _*(p(x,a))$. If $p(x,a)$
\th-forks over $A$ we get the conclusion of the claim by
definition of $\uth$-rank (after extending the types so we get
\th-division without changing the ranks). We will show that this
is the only possible case.

Suppose that $p(x,a)$ does not \th-fork over $A$. Let $p(x,a,B)$
be a non-\th-forking extension of $p(x,a)$ to $Ba$. Since
\th-forking is preserved by automorphisms, we know that for any
$a'\models tp(a/B)$ $p(x,a',B)$ is a non-\th-forking extension of
$p(x,a')$ and
\[ p(x,a',B)\upharpoonright B= p(x,a,B)\upharpoonright B :=p(x,B). \]
Given any such $a'$, $p(x,a',B)$ does not \th-fork over $A$ by
transitivity and by monotonicity $p(x,a',B)$ is a non-\th-forking
extension of $p(x,B)$.

Now, applying lemma \ref{amalgamation} $k$ times, we can find a
sequence $<a_1, a_2, \dots , a_k>$ with $a=a_1$ such that for any
$i\leq k$, $a_i\models tp(a/B)$ and for all $i,j, i\neq j$, $
a_i\thind _B a_j$ and $\cup_{i=1}^k p(x,a_i,B)$ is a
non-\th-forking extension of $p(x,B)$. But restricting to
$A\bar{a_i}$ we get that $\cup_{i=1}^k p(x,a_i)$ is a
non-\th-forking extension of $p(x)$ and therefore a
non-\th-forking extension of $p(x,a_1)$, a contradiction.

\end{proof}

\section{\th-forking in Simple and O-minimal Theories}

In this chapter we will look at theories for which independence
notions have been studied by other authors. In particular, we are
interested in analyzing the behavior of \th-forking inside simple
and o-minimal theories which are, as we mentioned in the
introduction, the two classes of theories where the independence
notion has been a main tool for the development of the respective
theories.

\subsection{Simple and Stable cases}

Given any theory, it is clear by the definition of \th-forking
that any formula which \th-divides over some set $A$ \th-forks
over $A$; Also by definition, for any type $p(x)$, any two sets of
formulae $\Delta$, $\Pi$ and any $k\in \mathbb{N}$,
\th$(p,\Delta,\Pi,k)\leq D(p,\Delta, k)$. So all simple theories
are rosy theories and \th-forking is a stronger notion than
forking in the sense that for any sets $A,B,C$, $A\ind_C
B\Rightarrow A\thind _C B$. As for the converse, we will prove in
this section that in all known examples of simple theories the two
notions agree. In fact we prove that the two notions being equal
is implied by the (strong) stable forking conjecture.

\subsubsection{\th-forking inside stable theories}

We will prove that in a stable theory \th-forking and forking are
the same. The following theorem will be the key to our proof.

\begin{theorem}
Inside a stable theory $T$, if a type $p$ forks over $A$, then
there is a non-forking extension $q$ of $p$, a formula $\delta $
and a non-algebraic type $\pi $ (in $\mathcal{L}^{eq}$) such that
$\left\{ \delta \left( x,a\right) \right\} _{a\models \pi }$ is
$k$ -inconsistent and $\delta \left( x,a\right) \in q$ for some
$a\models \pi $.
\end{theorem}

\begin{proof}
Let $\phi \left( x,b_{0}\right) $ be a formula implied by $p$
which divides over $A,$ and let $\left\{ \phi \left(
x,b_{i}\right) \right\} $ be the set of $k$-inconsistent formulas
which witness the dividing, with $\left< b_{i} \right>$ an $A$
indiscernible sequence. Since we are working inside a stable
theory, we can assume $k=2.$ Let $\Phi \left( x,y\right) =\left\{
\phi (x,y),\text{ }x=y\right\} ,$ $\pi_{1}=tp\left( b_{0}/A\right)
.$

Since we are only looking for a non-forking extension of $p$, we
can assume $p$ to be a complete type in $S(acl(A))$ so that
$\left( R,d\right) _{\Phi }\left( p\right) =\left( \alpha
,1\right) $. Let $\phi_0$ be a $\Phi $ formula in $p$ with the
same $\Phi $-rank and degree as $p$. Then $\phi _{0}\wedge \phi
\left( x,b_{0}\right) $ we be an $A$-dividing $\Phi$-formula such
that $(R,d)_\Phi (p)= (R,d)_\Phi (\phi_1)$; we can assume that
$(R,d)_{\Phi} (p)=(R,d)_{\Phi} (\phi(x,b_0))$.

To complete the proof we will need some sort of pre-\th-dividing
in a definable class of formulas in $p$ (or in a non-forking
extension). We suspect this may have been done in the theory of
canonical bases in stability theory. However, we were not able to
find a statement in the literature with exactly the properties we
needed so we will prove the following claim.

\begin{claim}
Let $p$ be a type in $S(Ab_0)$ such that $\left( R,d\right) _{\Phi
}\left( p\right) =\left( \alpha ,1\right) .$ Let $\phi \left(
x,b_{0}\right) \in p$ such that $\left( R,d\right) _{\Phi }\left(
\phi \left( x,b_{0}\right) \right) =\left( R,d\right) _{\Phi
}\left( p\right) .$ Then we can find formula $\eta \left( x,\left[
b\right] \right) $ and a type $ q$ such that
\begin{enumerate}
\item $\eta \left( x,\left[ b\right] \right) $ belongs to $p\left(
x,b_{0}\right)$, \item $R_{\Phi}\left( \eta \left( x,\left[
b\right] \right) \right) =\alpha $, \item For any $[a],[b]$ such
that $\left[ a\right] ,\left[ b\right] \models q$ and $[a]\neq
[b]$,
\[R_{\Phi }\left(
p(x)\cup \left( \eta ( x,[ b] ) \wedge \eta ( x,[ a] ) \right)
\right) <R_{\Phi }\left( p\right) .\]
\end{enumerate}
If, besides, there is an infinite sequence $\left\langle
b_{i}\right\rangle $ of elements satisfying $tp(b_0/A)$ such that
$\left\{ \phi \left( x,b_{i}\right) \right\} $ is
$2$-inconsistent, then $tp\left( \left[ b\right] /A\right) $ is
non-algebraic.
\end{claim}

\begin{proof}
As in (see ~\cite{Ziegler} 2.2), given two definable classes
$\mathbb{F}$ and $\mathbb{G}$, we define
\begin{eqnarray}
\mathbb{F} \subset _{\Phi ,\alpha } \mathbb{G} \text{   if   }
R_\Phi (\mathbb{F}\backslash \mathbb{G})<\alpha
\end{eqnarray}
and
\begin{eqnarray}
\mathbb{F} \cong _{\Phi ,\alpha } \mathbb{G} \text{   if   }
R_\Phi (\mathbb{F}\triangle \mathbb{G})<\alpha.
\end{eqnarray}
Given formulas $\psi_1(x), \psi_2(x)$, we will abbreviate
$\psi_1(\mathbb{C})\subset_{\Phi, \alpha} \psi_2(\mathbb{C})$ by
$\psi_1(x)\subset_{\Phi, \alpha} \psi_2(x)$.

Let $\phi '$ be the formula $\phi \left( x,y\right) $ with
variable $y$ and parameters in $x.$ Let $\pi(y)$ be $tp(b_0/A)$.
For any $b\models \pi ,$ let $\psi \left(b, y\right) $ be the
$\phi $ definition of $p\cup \left\{ \phi \left( x,b\right)
\right\} $ (which is a $\phi '$ -formula over $Ab$). Then $\models
\psi(b,c)$ if and only if $\phi \left( x,b\right) \subset _{\Phi
,\alpha }\phi \left( x,c\right) $, for any $c$.  By using an
automorphism between $b$ and any other $a\models \pi ,$ we know
that $\psi _{a}\left( y\right) =\psi \left( a,y\right) .$ Let
\begin{equation*}
E_{1}\left( z,y\right) :=\psi \left( z,y\right)
\end{equation*}
For any $a\models \pi ,$ the $\Phi $-degree of $\phi \left(
x,a\right) $ is 1 so
\begin{eqnarray*}
\psi \left( a,b\right) &\iff &\phi \left( x,b\right) \subset
_{\Phi ,\alpha
}\phi \left( x,a\right) \\
&\iff &R_{\Phi }\left( \phi \left( x,b\right) \wedge \lnot \phi
\left(
x,a\right) \right) <\alpha \\
&\iff &R_{\Phi }\left( \phi \left( x,b\right) \wedge \phi \left(
x,a\right)
\right) =\alpha \\
&\iff &R_{\Phi }\left( \lnot \phi \left( x,b\right) \wedge \phi
\left(
x,a\right) \right) <\alpha \\
&\iff &\psi \left( b,a\right) \\
&\iff &\phi \left( x,a\right) \cong _{\Phi ,\alpha }\phi \left(
x,b\right)
\end{eqnarray*}

So $\pi \left( z\right) \wedge \pi \left( y\right) \vdash ``E_{1}
\text{ is an equivalence relation }$'' and there is some finite
$\pi _{1}\subset \pi $ such that $\pi _{1}\left( z\right) \wedge
\pi _{1}\left( y\right) \wedge E_{1}\left( z,y\right) $ is an
equivalence relation $E_{2}.$ However, we may have more points
than we want in each equivalence class so we restrict it a bit
more by finding a finite $\pi _{2}\subset \pi $ such that for some
(and because $\pi$ is a complete type over $A$ for all) $b\models
\pi ,$
\[ \left( R,d\right) _{\phi '}\left( \pi
_{2} (y) \cup \left\{ E_2\left( b,y\right) \right\} \right)
=\left( R,d\right) _{\phi '}\left( \pi (y) \cup \left\{ E_2\left(
b,y\right) \right\} \right) .
\] Let
\begin{equation*}
E(z,y):=\pi _{2}\left( z\right) \wedge \pi _{2}\left( y\right)
\wedge E_{2}\left( z,y\right).
\end{equation*}
$E(z,y)$ is contained in $\pi \left( z\right) \wedge \pi \left(
y\right) \wedge E\left( z,y\right) $ so for any $a,b\models \pi $

\begin{eqnarray}
E\left( a,b\right) \iff & E_{1}\left( a,b\right) \iff & R_{\Phi
}\left( \lnot \phi \left( x,a\right) \wedge \phi \left( x,b\right)
\right) <\alpha \\
& \lefteqn{\text{and, by additivity of $R_\Phi$-rank in stable
theories,}}
&  \nonumber \\
\lnot E\left( a,b\right) \Rightarrow & \lnot E_{1}\left(
a,b\right) \Rightarrow & R_{\Phi }\left( \phi \left( x,a\right)
\wedge \phi \left( x,b\right) \right) <\alpha ;
\end{eqnarray}
let $\left[ b\right] $ be the $E$-equivalence class of $b.$

We will now work with generic points to make calculations easier.
Let $b$ satisfy a non-forking extension of $tp\left( b_{0}/A\left[
b_{0}\right] \right) $ such that
$b\mathop{\mathpalette\Ind{}}_{A\left[ b_{0}\right] }b_{0}.$ Let
$c\models p\left( x,b_{0}\right) $ such that $b\thind _{Ab_{0}}c.$

Let $\eta (x)$ be the $\phi '$ definition of $tp\left(
b/Ab_{0}\right) $ (and because $tp(b/Ab_0)$ is a non-forking
extension of $tp\left( b/A\left[ b_{0}\right] \right) $ all the
parameters of $\eta$ are in $acl(A[b_0])$ so $\eta =\eta \left(
x,\left[ b_{0}\right] \right)$ is over $acl(A\left[ b_{0}\right] )
).$ In particular, we know that $b\models tp\left( b_{0}/A\left[
b_{0}\right] \right) $, $c\models p\left( x,b_{0}\right) $ and
$\psi \left( b_{0},y\right) \in tp\left( b/Ab_{0}\right) .$ It
follows that $\eta \left( x,\left[ b_{0}\right] \right) \in p$
(see ~\cite{Pillay} lemma 2.8).

We can now prove that such $\eta(x,[b_0])$ satisfies condition 3.
Suppose not and let $[b]\neq [a]$ such that $R_\Phi \left(
p(x)\cup (\eta (x, [b])\wedge \eta (x,[a])) \right)=\alpha$. We
extend $p(x)$ to a non-forking type $p'(x)$ over a model $M$.

Let $a',b'$ be elements in $[a]$ and $[b]$ respectively such that
$a'\ind _[a] A$ and $b'\ind _[b] A$. Let $c\models p(x)\cup (\eta
(x, [b])\wedge \eta (x,[a]))$ such that $c\ind _{M} a'b'$; such
$c$ satisfies the $\phi'$ definition of $tp(a'/A[a])$ and
$tp(b'/A[b])$. By ~\cite{Pillay} 2.8 again and non-forking, both
$b',a'\models \phi(c,y)$. But $[b']=[b]\neq [a]=[a']$ so $b'\ncong
_{\Phi ,\alpha }a'.$ By definition,
\begin{equation*}
\left( R,d\right) _{\Phi }\left( \phi \left( x,b'\right) \wedge
\phi \left( x,a'\right) \right) <\left( \alpha ,1\right)
\end{equation*} which contradicts $c\ind _M a'b'$.

Finally, if $\left\langle b_{i}\right\rangle $ is an infinite
sequence of elements satisfying $\pi $ such that $\left\{ \phi
\left( x,b_{i}\right) \right\} $ is $2$-inconsistent, then each
$b_{i}$ in the sequence belongs to a different class $\left[
b_{i}\right] $ which proves the last statement of the claim.
\end{proof}

We can prove the theorem using lemma \ref{amalgamation} and a
careful induction on the ranks. However, this proof may not work
in the general context of simple theories with stable forking
conjecture and it does not provide a specific \th-forking formula
which is always nice to have. We give a different proof which,
even though it is more technical, provides an algorithm to
construct the \th-dividing formula and it is a good insight of
what happens, for example, in algebraically closed fields.

We know that $\eta \left( x,\left[ b_{0}\right] \right) $ in the
claim forks over $A$ so there must be some formula in the $\eta
$-type of $p\left( x,b_{0}\right) \upharpoonright acl\left(
A\left[ b_{0}\right] \right)$ which divides over $A$. We can take
this new formula instead of $\phi \left( x,b_{0}\right) $ in the
statement of the theorem and assume that $ \phi \left( x,b\right)
$ itself satisfies the conclusion of the claim or that $\eta
\left( x,\left[ b_{0}\right] \right) $ is a $\Phi $-formula.

Given any two $\beta ,d\in N,$ the statement ``There is an
infinite subfamily $I$ of $ b_{i}\models \pi $ such that $\left(
R,d\right) _{\Phi }\left\{ \phi \left( x,b_{i}\right)
\right\}_{b_i\in I} \geq \left( \beta ,d\right) $'' is type
definable. If $\left\{ \phi \left( x,b\right) \right\}_{b\models
tp(b_0/A)} $ is not $k$-inconsistent there is an infinite
subfamily which is consistent. Let $\left( \beta _{1},d_{1}\right)
$ be the maximal $\Phi$-rank and -degree that an infinite
subfamily of $\{ \phi(x,b) \}_{b\models tp(b_0/A)}$ can
achieve\footnote{So that for any infinite subfamily $I'$ of
$\left\{ b\mid b\models tp(b_o/A) \right\} ,$
\[ \left( R,d\right)_{\Phi }
\left( \bigcup_{b\in I'} \left\{ \phi \left( x,b\right) \right\}
\right) \leq \left( \beta _{1},d_{1}\right).\] Such a maximum
exists by compactness}. Let $B$ be an infinite set of elements in
$tp(b_0/A)$ such that \[R_\Phi \left( \left\{ \phi \left(
x,b'\right) \right\}_{b'\in B} \right)=(\beta_1,d_1).\] By using
an automorphism we can assume that this set includes $\phi \left(
x,b_{0}\right) $. Note that by the previous claim $\beta _{1}$
must be less than $\alpha $.

The rest of the proof will consist of narrowing down the set of
sentences we will consider, and maybe working with non-forking
extensions.

We will define a sequence of formulas $\left\{ \theta _{i}\left(
x\right) \right\}$ inductively. Let $\phi_1(x,b)=\phi(x,b)$, $q_1$
be $\{\phi (x,b) \}_{b\in B}$ so $\left( R,d\right) _{\Phi }\left(
q_1 \right) =\left( \beta _{1},d_{1}\right).$

By the local property of the rank, there is a finite conjunction
of $\phi \left( x,b_{i}\right)  \in q_{1}$ (which we will call $
\theta _{1}$) such that $\left( R,d\right) _{\Phi }\left( \theta
_{1}\right) =\left( R,d\right) _{\Phi }\left( q_{1}\right)$; so
for any $b_{i},$ $\phi _{1}\left( x,b_{i}\right) \wedge \theta
_{1}$ is a forking extension of $\phi _{1}\left( x,b_{i}\right) $.
By extension, $p\left( x,b_{i}\right) \cup \left\{ \lnot \theta
_{1}\right\} $ is a non-forking extension $p_{2}\left(
x,b_{i}^{2}\right) $ of $p\left( x,b_{i}\right) .$ Let $r_1(x)$ be
a non-forking extension of $\theta_1$ to a model $M$ and let
$\gamma (y)$ be its $\phi$-definition. By compactness we can
assume $B$ is very big (big enough to use Erdos-Rado), and by
Erdos-Rado we can also assume it contains some infinite
$M$-indiscernible subset $< b_j >_{j\in J}$ (as usual we can get
$b_{j_0}=b_{0}$). For all $j\in J$ let $\phi _{2}\left(
x,b_{j}^{2}\right) =\phi \left( x,b_{j}\right) \wedge \lnot \theta
_{1}\left( x\right) $ which is consistent (and in fact a
non-forking extension of $\phi(x,b_j)$) for all $j\in J.$

Recall that $q_1(x)$ had the same $\Phi$ rank and degree as
$\theta_1$ so that $b_{j_0}=b_0\models \gamma(y)$. We will say
that $b_i^2\models \gamma(y)$ whenever the projection $b_j\models
\gamma(y)$.

Let $\pi _{2}=tp\left( b_{0}^{2}/M\right) \cup \gamma(y) $. This
type is satisfied by all the $b_i^2$'s and is therefore non
algebraic.

For any infinite subfamily $s(x)$ of $\left\{ \phi _{2}\left(
x,b_{i}^{2}\right) \right\} _{b_{i}^{2}\models \pi _{2}},$ we know
that $s(x)$ is of the form $\left\{ \phi \left( x,b_{i}\right)
\right\}_{i\in J} \cup \{ \lnot \theta _{1}(x) \} $ where $J$ is
some index set such that $b_i\models tp(b_0/M)$ for any $i\in J$.
By maximality of $(R,d)_{\Phi}(q_{1}) $ we have $\left( R,d\right)
_{\Phi }\left( r\right) \leq \left( \beta _{1},d_{1}\right) $. If
$ \left\{ \phi _{2}\left( x,b_{i}^{2}\right) \right\}
_{b_{i}^{2}\models \pi _{2}}$ is not $k$-inconsistent for any $k$,
we can once more find some consistent subset $I_2\subset
\pi_2(\mathcal{C})$ such that $\left( \bigwedge _{i\in I_2}
\phi_2(x,b_i^2) \right)$ has  maximum $\left( R,d\right)_\Phi$ in
the lexicographical order; say $\left( R,d\right) \left( \bigwedge
_{i\in I_2} \phi_2(x,b_i^2) \right)=\left( \beta _{2},d_{2}\right)
\leq \left( \beta _{1},d_{1}\right) $.

Recall that for $i\in I_2$, $b_i^2\models \gamma(y)$ (the $\phi$
definition of $\theta_1$) so
\[ (R,d)_\Phi \left( \left\{ \phi (x,b_i) \right\} _{i\in I_2} \cup
\theta_1(x) \right)=(\beta_1, d_3)\] for some $d_3.$ If
$\beta_1=\beta_2$ we get, using the fact that $\theta_1$ is a
$\Phi$-formula, that

\begin{eqnarray*}
\deg_{\Phi} \left( \left\{ \phi \left( x,b_{i}\right)
\right\}_{i\in I_2} \right) && =\deg _{\Phi }\left( \left\{ \phi_2
\left( x,b^2_{i}\right) \right\}_{i\in I_2} \right) +\deg_\Phi
\left( \left\{ \phi \left( x,b_{i}\right)
\right\}_{i\in I_2} \cup \theta_1(x) \right)  \\
&& \Rightarrow d_{1}\geq d_{2}+d_3 \\
&& \Rightarrow d_2<d_1. \\
\end{eqnarray*} So $(\beta_2, d_2) <^{lex} (\beta_1, d_1)$.
Let $q_{2}=\left\{ \phi (x,b_i^2) \right\}_{i\in I_2}$ and $\theta
_{2}$ be a conjunction of $\phi _{2}\left( x,b_{i}^{2}\right) $
with the same $\Phi$-rank, degree as $q_{2}$.

We can repeat the above construction as long as \[ \left\{ \phi
_{n}\left( x,b_{i}^{n}\right) \right\} _{b_{i}^{n}\models \pi
_{n}}\] is not $k$-inconsistent for any $k$; by doing this we
would get a sequence $\left< (R,d)_\Phi \left( \theta_i \right)
\right>_{i\leq n}$ which is strictly decreasing in the
lexicographical order. Since both the $\Phi$ rank and degrees are
ordinals (and therefore well ordered) at some point we must get
some $N$ such that $\left\{ \phi _{N}\left( x,b_{i}^{N}\right)
\right\} _{b_{i}^{N}\models \pi _{N}}$ is $k$-inconsistent and by
construction, such $\phi _{N}\left( x,b_{0}^{N}\right)$ would
belong to a non-forking extension of $p\left( x,b_{0}\right) .$
\end{proof}

\begin{corollarytotheproof}
\emph{If in some model of a simple theory, some type $p$ has a
stable formula which forks over $A$, then there is a non-forking
extension $q$ of $p$, a stable formula $\delta $ and a type $\pi $
such that $\left\{ \delta \left( x,a\right) \right\} _{a\models
\pi }$ is $k$-inconsistent and $\delta \left( x,a\right) \in q$
for some $a\models \pi $}
\end{corollarytotheproof}

\begin{proof}
If $\phi$ is stable, then so is the formula $\eta$ defined in the
claim (it is equivalent to a conjunction of $\delta$ formulas) and
so are all the formulas defined in the rest of the proof.
\end{proof}

\begin{definition}
A theory $T$ has the \emph{strong stable forking property} if
whenever $p$ is a type over some $B$ subset of $\mathcal{C}$ and
$p$ forks over a set $A$ (not necessarily contained in $B$), then
this is witnessed by some instance of a stable formula in $p$. We
will say that $T$ satisfies the \emph{weak stable forking
property} if we require $B$ to be a model of $T$ and $A$ to be a
subset of $B$.
\end{definition}

\begin{corollary}
\label{stableforkingconjecture} In any stable theory $T$, if
$M\models T$ is a model of $T$, $A$ is a subset of $M$ , $q\in
S(M)$ is a complete type and $p=q\upharpoonright A$, then $q$ is a
non-forking extension of $p$ if and only if it is a
non-\th-forking extension. In fact, this conclusion is true for
any simple theory $T$ for which the strong stable forking property
holds \footnote{Even though the proof is done assuming the strong
stable forking property, we think the weak stable forking
conjecture might be enough.} (which includes all known examples of
simple theories).
\end{corollary}

\begin{proof}
The left to right implication is immediate since \th-forking is a
stronger notion than forking. For the other side, if $q$ forked
over $A$ there would be some $r\left( x,b'\right) $ a non-forking
extension of $q$, a formula $\delta $ and a non-algebraic type
$\pi $ such that $\left\{ \delta \left( x,b\right) \right\}
_{b_{i}\models \pi }$ is $k$-inconsistent, $b'\models \pi ,$
$\delta \left( x,b'\right) \in r\left( x,b'\right) $. Let $\left\{
b_{i}\right\} _{i\in N}$ be an infinite set of $A$-conjugates of
$b'$. By monotonicity,
\[ \text{\th} \left( p\left( x\right) \wedge \delta \left( x,b_{i}\right)
,\delta ,\pi ,k\right) \geq \text{\th} \left( r\left(
x,b_{i}\right) ,\delta ,\pi ,k\right). \] Conjugating over $A$ we
get \th $\left( r\left( x,b' \right) ,\delta , \pi ,k\right) =$\th
$\left( r\left( x,b_{i}\right) ,\delta , \pi ,k\right) $. By
definition,
\[
\text{\th} \left( p,\delta ,\pi ,k\right) \geq \text{\th} \left(
r\left( x,b' \right) ,\delta ,\pi ,k\right) +1.\] But any
non-forking extension is a non-\th-forking extension so
\begin{equation*}
\text{\th} \left( r\left( x,b' \right) ,\delta ,\pi ,k\right) =
\text{\th }\left( q,\delta , \pi ,k\right)
\end{equation*} so
\[
\text{\th} \left( p ,\delta ,\pi ,k\right) < \text{\th }\left(
q,\delta , \pi ,k\right). \] \end{proof}

\subsubsection{Supersimple theories}

In this section we shall prove that for supersimple theories the
forking relation is equivalent to the \th-forking one. It is clear
that for any type $p$, $\text{U}(p)\geq \uth (p)$; so in
supersimple theories both the $\text{U}$ and the $\uth$ are
defined and they both characterize forking and \th-forking
respectively. Therefore, it is enough to prove the inequality
assuming both these ranks are ordinal valued.

\begin{theorem}
For any type $p$
\[ \text{U}(p)=\uth (p) \]
\end{theorem}

\begin{proof}
The U-rank is always bigger than the $\uth$-rank and we proved the
equivalence between $\uth _* $ and $\uth$, so it is enough to show
that for any ordinal $\alpha$ and any type $p$,
\[
\text{U}(p)=\alpha \Rightarrow \uth _* (p)\geq \alpha.
\]

We will do an induction on U$(p)$. For $\alpha=0$ or a limit
ordinal it is clear. Let us assume then that
U$(p)=\alpha=\sigma+1$ and that for any $\lambda \leq \sigma$ and
any type $p$, $\text{U}(p)=\lambda \Rightarrow \uth _* (p)\geq
\lambda$. But $\text{U}(p)\geq \uth_*(p)$ so our induction
hypothesis is equivalent to $\text{U}(p)=\lambda \Rightarrow \uth
_* (p)= \lambda$.

Let $p$ be a type over some set $A$ such that $\text{U}(p)=
\sigma+1$. By definition, we can extend $p$ to some type $q(x,C)$
over $C\supset A$ such that $\text{U}(q)= \sigma$ and $q$ forks
over $A$. Let $a'$ be the canonical base of $q$ over
$A$\footnote{The construction of such hyperimaginary would be made
by defining an equivalence relation just for extensions of $p(x)$
(looking at $A$ as a fixed set of parameters). We might have to
extend $p(x)$ to a type over a model containing $A$, but all this
can be done by the extension and transitivity properties.} and let
$p(x,a')=q\upharpoonleft Aa'$. Now, by elimination of
hyperimaginaries in supersimple theories
(see~\cite{Buechler/Pillay/Wagner}), we can assume $a'$ is a tuple
in the model; by definition of the canonical base we know that for
any two $a_1,a_2\models tp(a'/A)$, $p(x,a_1)\cup p(x,a_2)$ is a
forking extension of $p(x,a_1)$ and by \cite{Kim} theorem 3.3.6 we
have that $q(x)$ does not fork over $Aa'$. Now, in supersimple
theories this translates into the following statement. For any
$a_1,a_2\models tp(a'/A), $
\[ \sigma=\text{U}(q)=\text{U}\left( p(x,a_1)\right) >\text{U}\left(
p(x,a_1)\cup p(x,a_2) \right).
\]
By induction hypothesis,
\[ \sigma=\uth _* \left( p(x,a_1)\right) >\uth _* \left( p(x,a_1)\cup
p(x,a_2) \right). \] By definition of the $\uth _*$ we get $\uth
_* (p) \geq \sigma +1=\text{U}(p)\geq \uth(p)=\uth_*(p)$. So
$\uth_*(p)=\sigma +1=\text{U}(p)$.
\end{proof}

\subsection{O-minimal:}

We shall prove that in o-minimal theories the global $\uth$ rank
is defined and equal to the usual dimension. This will prove that
$\thind$ is precisely the usual independence relation. Let $M$ be
the model of some o-minimal theory.

 \begin{claim}
 For any definable $A\subset M^k$, if $\dim (A)=n$ then $\uth (A)\leq n$
 \end{claim}

 \begin{proof}

Suppose otherwise and let $i$ be the smallest integer such that
there is a definable $A$, $\dim (A)=i$ and $\uth (A)\geq i+1$, let
$A$ be definable over some finite set $B$ by some formula $\phi$,
and let $A\subseteq M^k$. We can then find some ``generic'' point
(in the sense of the $\uth$-rank) $a\in A$ such that $\uth \left(
tp(a/B) \right) \geq i+1$. By minimality of $i$, $\dim
(tp(a/B))\geq i$ (otherwise we could find some definable subset
contained in $A$ for which $a$ is dimension-generic). By the
theory of o-minimal sets (see \cite{vandenDries}) there is some
open subset $A'$ with $a\in A'\subset A$ such that we can find a
projection from $M^k$ to $M^i$ which is a definable homomorphism
from $A'$ into some open subset of $M^i$. But we know that $\uth
(A')\geq i+1$ and $\dim (A')=i$, both of which are preserved by
definable homomorphisms so we may assume that $A\subseteq M^i$.

By definition of the $\uth$-rank, there is some formula $\delta
(x,y)$, an element $b_0$ and some $C\supset B$ such that
\begin{itemize}
\item $tp(b_0/C)$ is non algebraic \item $\uth \left( tp(a/B)\cup
\delta (x,b_0)\right) \geq i$ and \item $\{ \delta (x,b')
\}_{b'\models tp(b_0/C)}$ is $k$-inconsistent.
\end{itemize}

We will first do the case where $i=1$. Let $\left< b_i \right>
_{i\in \mathbb{N}}$ be some indiscernible sequence in $tp(b_o/C)$
and $A_j$ be the set defined by $\delta (x,b_j)$. By hypothesis
$\uth(A_j)=1$ and $A_j\subset M^1$ so all of the $A_j$'s are
infinite subsets of an o-minimal model. By definition of
o-minimality, each of the $A_j$'s must contain an interval. We
could therefore (uniformly) define $A_j'$ as the first interval
contained in $A_j$ and we can define the subset $X\subset M$
consisting of left endpoints of the $A_j'$. By $k$-inconsistency,
this set cannot contain an interval and it must be infinite (not
more than $k-1$ of the $A_j$'s can have the same left endpoint)
contradicting o-minimality.

For a general $i$, note that by minimality of $i$ we have that
$\dim (\phi(x)\wedge \delta(x,b'))\geq i$ for any $b'\models
tp(b/C)$. But it is a definable subset of $M^i$ so the projection
to each of the coordinates must be infinite and definable. Taking
the projection to the first coordinate
\[\phi_0(x_0,b'):=\exists x \phi(x)\wedge \delta (x,b')\wedge
\pi_0(x)=x_0\] we have a definable set in $M^1$ so by the $i=1$
case the set of formulas $\{ \phi_0 (x_0, b' )\}_{b'\models
tp(b_0/C)}$ cannot be $k$-inconsistent. By definition of
$k$-inconsistency we can find an infinite subsequence $\left< b_j
\right> _{j\in J}$ of $\{ b' \}_{b'\models tp(b_0/B)}$ such that
\[\bigwedge_{b_i} \phi_0(x_0, b_i) \] is consistent; say it is
realized by $a_0$. Let $q_0(y)$ be the (non algebraic) type
$tp(b_0/B)\cup \{\phi_0(a_0,y)\}$ and consider $\{ \delta
(x,b')\}_{b'\models q_0}$. We repeat the procedure with each of
the projections (we always get non algebraic types) and
inductively define a point in $M^i$ $a:=(a_0, \dots , a_{i-1})$
which realizes infinitely many of the $\delta (x,b')$'s,
contradicting $k$-inconsistency.
\end{proof}

The other inequality is much easier. If $A\subset M^k$ and $\dim
(A)=n$ then there is some open subset $A'$ of $A$ and a definable
projection $\pi$ from $M^k$ to $M^n$ such that $A'$ gets sent to
an open subset $U$ of $M^n$ and $\pi :A'\rightarrow U$ is a
homomorphism. Thus, $\dim (U)=n$. This implies that the projection
to each of its coordinates is infinite, but $x_i=a$ is clearly a
\th-forking formula for any $i\leq n$ which proves $\uth (U)\geq
n$ so $\uth (A')\geq n$. By monotonicity of the rank, $\uth
(A)\geq n$.





\nocite{vandenDries, Kim/Pillay/stableforking, hodgess}

\bibliographystyle{alpha}

\bibliography{jsl}

\begin{thebibliography}{BPW01}

\bibitem[BPW01]{Buechler/Pillay/Wagner}
Steven Buechler, Anand Pillay, and Frank Wagner.
\newblock Supersimple theories.
\newblock {\em J. Amer. Math. Soc.}, 14(1):109--124 (electronic), 2001.

\bibitem[CW02]{Casanovas/Wagner}
Enrique Casanovas and Frank Wagner.
\newblock The free roots of the complete graph.
\newblock {\em Submitted for publication}, 2002.

\bibitem[Hod93]{hodges}
Wilfrid Hodges.
\newblock {\em Model theory}, volume~42 of {\em Encyclopedia of Mathematics and
  its Applications}.
\newblock Cambridge University Press, Cambridge, 1993.

\bibitem[Hod97]{hodgess}
Wilfrid Hodges.
\newblock {\em A shorter model theory}.
\newblock Cambridge University Press, Cambridge, 1997.

\bibitem[Kim96]{Kim}
Byunghan Kim.
\newblock {\em Simple First Order Theories}.
\newblock PhD thesis, University of Notre Damme, 1996.

\bibitem[KP97]{Kim/Pillay}
Byunghan Kim and Anand Pillay.
\newblock Simple theories.
\newblock {\em Ann. Pure Appl. Logic}, 88(2-3):149--164, 1997.
\newblock Joint AILA-KGS Model Theory Meeting (Florence, 1995).

\bibitem[KP01]{Kim/Pillay/stableforking}
Byunghan Kim and A.~Pillay.
\newblock Around stable forking.
\newblock {\em Fund. Math.}, 170(1-2):107--118, 2001.
\newblock Dedicated to the memory of Jerzy \L o\'s.

\bibitem[Pil96]{Pillay}
Anand Pillay.
\newblock {\em Geometric stability theory}, volume~32 of {\em Oxford Logic
  Guides}.
\newblock The Clarendon Press Oxford University Press, New York, 1996.
\newblock Oxford Science Publications.

\bibitem[She90]{Shelah}
S.~Shelah.
\newblock {\em Classification theory and the number of nonisomorphic models},
  volume~92 of {\em Studies in Logic and the Foundations of Mathematics}.
\newblock North-Holland Publishing Co., Amsterdam, second edition, 1990.

\bibitem[vdD98]{vandenDries}
Lou van~den Dries.
\newblock {\em Tame topology and o-minimal structures}, volume 248 of {\em
  London Mathematical Society Lecture Note Series}.
\newblock Cambridge University Press, Cambridge, 1998.

\bibitem[Wag00]{Wagner}
Frank~O. Wagner.
\newblock {\em Simple theories}, volume 503 of {\em Mathematics and its
  Applications}.
\newblock Kluwer Academic Publishers, Dordrecht, 2000.

\bibitem[Zie98]{Ziegler}
Martin Ziegler.
\newblock Introduction to stability theory and {M}orley rank.
\newblock In {\em Model theory and algebraic geometry}, volume 1696 of {\em
  Lecture Notes in Math.}, pages 19--44. Springer, Berlin, 1998.

\end{thebibliography}


\end{document}